\documentclass[11pt, reqno]{amsart}

\usepackage[OT2, T1]{fontenc}

\usepackage[usenames,dvipsnames]{color}

\usepackage{url}
\usepackage{amsmath}
\usepackage{array}
\usepackage{graphicx}
\usepackage{amsfonts}
\usepackage{amssymb}
\usepackage{amsthm}
\usepackage{hyperref}
\usepackage{paralist}
\usepackage{colonequals}	
\usepackage{color}
\usepackage{mathabx}		
\usepackage{amscd}
\usepackage[all,cmtip]{xy}	
\usepackage{sseq}			
\usepackage{verbatim}
\usepackage{parskip}
\usepackage{microtype}
\usepackage[in]{fullpage}
\usepackage{mathrsfs} 			
\usepackage{bm} 				
\usepackage{cleveref}
\usepackage{enumitem}
\usepackage{moreenum}			
\usepackage[alphabetic,lite]{amsrefs} 	
\usepackage{stmaryrd}	

\DeclareSymbolFont{cyrletters}{OT2}{wncyr}{m}{n}
\DeclareMathSymbol{\Sha}{\mathalpha}{cyrletters}{"58}

\newcommand{\gA}{\alpha}

\newcommand{\bA}{\mathbb{A}}
\newcommand{\bC}{\mathbb{C}}

\newcommand{\bF}{\mathbb{F}}
\newcommand{\bG}{\mathbb{G}}

\newcommand{\bP}{\mathbb{P}}

\newcommand{\bR}{\mathbb{R}}

\newcommand{\bZ}{\mathbb{Z}}

\newcommand{\bbB}{\mathbf{B}}

\newcommand{\bbR}{\mathbf{R}}

\newcommand{\cA}{\mathcal{A}}

\newcommand{\cC}{\mathcal{C}}

\newcommand{\cF}{\mathcal{F}}
\newcommand{\cG}{\mathcal{G}}
\newcommand{\cH}{\mathcal{H}}

\newcommand{\cO}{\mathcal{O}}
\newcommand{\cP}{\mathcal{P}}

\newcommand{\fa}{\mathfrak{a}}

\newcommand{\fm}{\mathfrak{m}}
\newcommand{\fo}{\mathfrak{o}}

\newcommand{\fS}{\mathfrak{S}}

\newcommand{\sF}{\mathscr{F}}

\newcommand{\sX}{\mathscr{X}}
\newcommand{\sY}{\mathscr{Y}}
\newcommand{\sZ}{\mathscr{Z}}

\newcommand{\ra}{\rightarrow}

\newcommand{\xra}{\xrightarrow}

\newcommand{\hra}{\hookrightarrow}

\newcommand{\wt}{\widetilde}
\newcommand{\wh}{\widehat}

\newcommand{\pr}{^{\prime}}

\newcommand{\ce}{\colonequals}
\newcommand{\ov}{\overline}

\renewcommand{\b}{\textbf}
\newcommand{\surjects}{\twoheadrightarrow}

\newcommand{\union}{\cup} 			
\newcommand{\tensor}{\otimes} 		
\newcommand{\isomto}{\overset{\sim}{\longrightarrow}}

\newcommand{\fppf}{\mathrm{fppf}}		
\newcommand{\et}{\mathrm{\acute{e}t}}	

\renewcommand{\i}{^{-1}}
\renewcommand{\th}{^{\mathrm{th}}}

\providecommand{\f}[2]{\frac{#1}{#2}}
\DeclareMathOperator{\Ker}{Ker}			
\DeclareMathOperator{\im}{Im}			
\DeclareMathOperator{\Spec}{Spec}		
\DeclareMathOperator{\Hom}{Hom}			
\DeclareMathOperator{\Char}{char}		
\DeclareMathOperator{\Frac}{Frac}		
\DeclareMathOperator{\id}{id}			
\DeclareMathOperator{\Ob}{Ob}			
\DeclareMathOperator{\Res}{Res}		
\DeclareMathOperator{\GL}{GL}		
\DeclareMathOperator{\Aut}{Aut}		
\DeclareMathOperator{\Isom}{Isom}		
\DeclareMathOperator{\SEC}{SEC}		
\DeclareMathOperator{\ET}{ET}		
\newcommand{\ba}{\begin{aligned}}
\newcommand{\ea}{\end{aligned}}
\newcommand{\be}{\begin{equation}}
\newcommand{\ee}{\end{equation}}
\newcommand{\pf}{\begin{proof}}
\newcommand{\bpf}{\begin{proof}}
\newcommand{\epf}{\end{proof}}
\newcommand{\bthm}{\begin{thm}}
\newcommand{\ethm}{\end{thm}}
\newcommand{\bprop}{\begin{prop}}
\newcommand{\eprop}{\end{prop}}
\newcommand{\bcor}{\begin{cor}}
\newcommand{\ecor}{\end{cor}}
\newcommand{\brem}{\begin{rem}}
\newcommand{\erem}{\end{rem}}
\newcommand{\brems}{\begin{rems} \hfill \begin{enumerate}[label=\b{\thesubsection.},ref=\thesubsection]}
\newcommand{\remi}{\addtocounter{subsection}{1} \item}
\newcommand{\erems}{\end{enumerate} \end{rems}}
\newcommand{\blem}{\begin{lemma}}
\newcommand{\elem}{\end{lemma}}
\newcommand{\bconj}{\begin{conj}}
\newcommand{\econj}{\end{conj}}
\newcommand{\bprob}{\begin{Problem}}
\newcommand{\eprob}{\end{Problem}}
\newcommand{\bq}{\begin{q}}
\newcommand{\eq}{\end{q}}
\newcommand{\benum}{\begin{enumerate}[label={(\alph*)}]}
\newcommand{\benuma}{\begin{enumerate}[label={(\arabic*)}]}

\newcommand{\eenum}{\end{enumerate}}
\newcommand{\bc}{}
\newcommand{\beg}{\begin{eg}}
\newcommand{\eeg}{\end{eg}}
\newcommand{\bcl}{\begin{claim}}
\newcommand{\ecl}{\end{claim}}
\newcommand{\lab}{\label}
\newcommand{\q}{\quad}
\newcommand{\qq}{\quad\quad}
\newcommand{\qqq}{\quad\quad\quad}
\newcommand{\qqqq}{\quad\quad\quad\quad}

\newcommand*{\QED}{\hfill\ensuremath{\qed}}
\newcommand*{\QEDD}{\hfill\ensuremath{\qed\qed}}

\theoremstyle{plain}
\newtheorem{thm}[subsection]{Theorem}
\Crefname{thm}{Theorem}{Theorems}

\Crefname{rethm}{Theorem}{Theorem}
\newtheorem{prop}[subsection]{Proposition}
\Crefname{prop}{Proposition}{Propositions}
\Crefname{eg}{Example}{Examples}
\newtheorem{Problem}[subsection]{Problem}
\Crefname{Problem}{Problem}{Problems}
\newtheorem{conj}[subsection]{Conjecture}
\Crefname{conj}{Conjecture}{Conjectures}
\newtheorem{cor}[subsection]{Corollary}
\Crefname{cor}{Corollary}{Corollaries}
\newtheorem{lem}[equation]{Lemma}
\Crefname{lem}{Lemma}{Lemmas}

\newtheorem{lemma}[subsection]{Lemma}

\newtheorem{sublem}[equation]{Lemma}
\Crefname{subprop}{Proposition}{Propositions}
\Crefname{subcor}{Corollary}{Corollaries}
\Crefname{sublem}{Lemma}{Lemmas}

\theoremstyle{remark}
\newtheorem{claim}[equation]{Claim}
\Crefname{claim}{Claim}{Claims}
\newtheorem{subrem}[equation]{Remark}
\Crefname{subrem}{Remark}{Remarks}

\theoremstyle{definition}
\newtheorem{eg}[subsection]{Example}
\newtheorem{rem}[subsection]{Remark}
\Crefname{rem}{Remark}{Remarks}
\newtheorem*{rems}{Remarks}

\newtheoremstyle{subsection-tweak}
   {11pt}
   {3pt}%
   {}
   {}%
   {\bfseries}
   {}%
   {.5em}
   {\thmnumber{\@{#1}{}\@{#2}.}%
    \thmnote{~{\bfseries#3.}}}

\Crefname{innercustomconj}{Conjecture}{Conjecture}

\theoremstyle{subsection-tweak}
\newtheorem{pp}[subsection]{}
\newcommand{\bpp}{\begin{pp}}
\newcommand{\epp}{\end{pp}}

\numberwithin{equation}{subsection}

\begin{document}
\author{K\k{e}stutis \v{C}esnavi\v{c}ius}
\title{Topology on cohomology of local fields}
\date{\today}
\address{Department of Mathematics, University of California, Berkeley, CA 94720-3840, USA}
\email{kestutis@berkeley.edu}
\urladdr{http://math.berkeley.edu/~kestutis/}


\begin{abstract} 
Arithmetic duality theorems over a local field $k$ are delicate to prove if $\Char k > 0$. In this case, the proofs often exploit topologies carried by the cohomology groups $H^n(k, G)$ for commutative finite type $k$-group schemes $G$. These ``\v{C}ech topologies'', defined using \v{C}ech cohomology, are impractical due to the lack of proofs of their basic properties, such as continuity of connecting maps in long exact sequences. We propose another way to topologize $H^n(k, G)$: in the key case $n = 1$, identify $H^1(k, G)$ with the set of isomorphism classes of objects of the groupoid of $k$-points of the classifying stack $\bbB G$ and invoke Moret-Bailly's general method of topologizing $k$-points of locally of finite type $k$-algebraic stacks. Geometric arguments prove that these ``classifying stack topologies'' enjoy the properties expected from the \v{C}ech topologies. With this as the key input, we prove that the \v{C}ech and the classifying stack topologies actually agree. The expected properties of the \v{C}ech topologies follow, which streamlines a number of arithmetic duality proofs given elsewhere.
\end{abstract}


\maketitle

\section{Introduction}

\bpp[The need of topology on cohomology]
Let $k$ be a nonarchimedean local field of characteristic $p \ge 0$. The study of cohomology groups $H^n(k, G)$ for commutative finite $k$-group schemes $G$ is facilitated by Tate local duality (extended to the $p \mid \#G$ case in \cite{Sha64}): 
\be\lab{TLD} \tag{$\ddag$}
H^n(k, G)\quad \text{and} \quad H^{2 - n}(k, G^D)\quad \text{ are Pontryagin duals, where $G^D$ is the Cartier dual of $G$.}
\ee
If $p \nmid \#G$, then $H^n(k, G)$ is finite, but if $p \mid \#G$ (when $H^n$ abbreviates $H^n_\fppf$), then such finiteness fails: for instance,
\[
 H^1(\bF_p((t)), \alpha_p) \cong \bF_p((t))/\bF_p((t))^p.
\]
Therefore, in general one endows $H^n(k, G)$ with a topology and interprets \eqref{TLD} as Pontryagin duality of Hausdorff locally compact abelian topological~groups.

The above illustrates a general phenomenon: local arithmetic duality theorems face complications in positive characteristic due to the failure of finiteness of various cohomology groups. These complications are overcome by exploiting topologies carried by the relevant cohomology groups. The aim of this paper is to present a new way to define these topologies. This way is outlined in \S\ref{our-app}~and
\begin{enumerate}[label={(\arabic*)}] 
\item \lab{a-priori}
\emph{A priori} seems more robust than the ``\v{C}ech cohomology approach'' used elsewhere;

\item \lab{a-posteriori}
\emph{A posteriori} gives the same topology as the usual ``\v{C}ech cohomology approach''.
\eenum
The interplay of \ref{a-priori} and \ref{a-posteriori} has particularly pleasant implications for local arithmetic duality theorems: as we explain in \S\ref{Mil-app}, it streamlines a number of proofs given elsewhere.
\epp

\bpp[The \v{C}ech cohomology approach] \lab{Mil-app}
This way to topologize $H^n(k, G)$, where $k$ is a nonarchimedean local field of characteristic $p > 0$ and $G$ is a commutative finite type $k$-group scheme, is explained in \cite{Mil06}*{III.\S6}. It first establishes a connection with \v{C}ech cohomology:
\[
\textstyle{H^n(k, G) = \varinjlim_L H^n(L/k, G)},\quad \text{where $L/k$ runs over finite extensions in an algebraic closure $\ov{k}$,}
\]
and then endows $H^n(L/k, G)$ with the subquotient topology of $G(\tensor_{i = 0}^n L)$, and $H^n(k, G)$ with the direct limit topology. This ``\v{C}ech topology'' on $H^n(k, G)$ is discussed further in \S\ref{Cech-topo}. 

A quotient topology or a direct limit topology can be difficult to work with individually, and the \v{C}ech topology, which combines them, seems particularly unwieldy. A number of its expected properties lack proofs, e.g., it is not clear whether $H^n(k, G)$ is a topological group, nor whether $H^n(k, G)$ is discrete if $G$ is smooth and $n \ge 1$, nor whether connecting maps in long exact cohomology sequences are continuous. Such continuity is used in the proofs of local arithmetic duality theorems in positive characteristic, see, for instance, \cite{Mil06}*{III.6.11 or III.7.8}. Not knowing it in general, one resorts to ad hoc modifications of these proofs to ensure their completeness.

In contrast, the way to topologize $H^n(k, G)$ outlined in \S\ref{our-app} permits geometric proofs of all the desired properties, which illustrates \ref{a-priori}. The corresponding properties of the \v{C}ech topology follow thanks to \ref{a-posteriori}. This renders the ad hoc modifications mentioned earlier unnecessary.
\epp

\bpp[The classifying stack approach] \lab{our-app}
This way to topologize $H^n(k, G)$, where $k$ and $G$ are as in \S\ref{Mil-app}, uses the identification of $H^1(k, G)$ with the set of isomorphism classes of $G$-torsors, i.e.,~the set of isomorphism classes of objects of the groupoid of $k$-points of the classifying stack $\bbB G$. In \cite{MB01}*{\S2}, Moret-Bailly showed how to exploit the topologies on $X(k)$ for locally of finite type $k$-schemes $X$ to topologize the set of isomorphism classes of objects of $\sX(k)$ for any locally of finite type $k$-algebraic stack $\sX$ (see \Cref{app-b}, esp.~\S\ref{stack-case}, for details). His definition applies to $\sX = \bbB G$ and hence topologizes $H^1(k, G)$. As for the other $H^n$, the identification $H^0(k, G) = G(k)$ should be the way to topologize $H^0$, whereas the topology on $H^n(k, G)$ for $n \ge 2$ could simply be defined to be discrete. The resulting ``classifying stack topology'' on $H^n(k, G)$ seems convenient to work with:
\benum
\item \lab{easy-a}
Its properties follow from geometric arguments, e.g.,~connecting maps from $H^0$ to $H^1$ are continuous because they have underlying morphisms of algebraic stacks; 

\item \lab{easy-b}
Its definition permits other local topological rings (e.g.,~the ring of integers of $k$) as bases, which allows us to treat all local fields simultaneously;

\item \lab{easy-c}
Its definition permits $G$ that are locally of finite type, or noncommutative, or algebraic spaces, or all of these at once.
\eenum
We exploit the properties of the classifying stack topology to prove in \Cref{H0-homeo} and \Cref{main-Cech,Hi-disc} that this topology actually agrees with the \v{C}ech topology. In particular, the latter is discrete for $n \ge 2$; this seems not to have been noticed in the existing literature.
\epp

\bpp[A summary of our conclusions] \lab{summary}
Due to \ref{easy-a}--\ref{easy-c}, the classifying stack topology lends itself to general settings. For convenience of a reader not interested in generalities (which occupy \S\S\ref{topo-coho}--\ref{cont-les}), we summarize our findings in the case of commutative finite type group schemes $G$ over local fields~$k$.

The $H^n(k, G)$ are locally compact Hausdorff abelian topological groups that are discrete for $n \ge 2$. If $G$ is smooth (in particular, if $\Char k = 0$), then $H^n(k, G)$ is discrete for $n \ge 1$. If $k$ is nonarchimedean, $\fo$ is its ring of integers, and $\cG$ is a finite type flat $\fo$-model of $G$, then $H^1(\fo, \cG)$ is compact and $H^1(\fo, \cG) \ra H^1(k, G)$ is continuous and open; if $\cG$ is separated, then $H^1(\fo, \cG)$ is Hausdorff.

For a short exact sequence 
\[
0 \ra H \ra G \ra Q \ra 0
\]
of commutative finite type $k$-group schemes, all the maps in its long exact cohomology sequence are continuous. 
\begin{itemize}
\item 
If $H$ is smooth (resp.,~proper), then $G(k) \ra Q(k)$ is open (resp.,~closed). 

\item
If $G$ is smooth (resp.,~proper), then $Q(k) \ra H^1(k, H)$ is open (resp.,~closed). 

\item
If $Q$ is smooth (resp.,~proper), then $H^1(k, H) \ra H^1(k, G)$ is open (resp.,~closed). 

\item
Finally, $H^1(k, G) \ra H^1(k, Q)$ is always open (but possibly not closed). 
\end{itemize}
All the subsequent maps in the long exact cohomology sequence are both open and closed.

Thanks to \Cref{main-Cech,Hi-disc}, all the claims above, except possibly the ones involving $\fo$, also hold for the \v{C}ech topology, for which a number of them seem not to have appeared in the literature.
\epp

\bpp[\'{E}tale-openness, proper-closedness, and finite-closedness] \lab{et-op-ann}
To take advantage of the flexibility provided by \S\ref{our-app} \ref{easy-a}--\ref{easy-b}, e.g., to be able to treat all local fields without distinguishing archimedean cases, we adopt the axiomatic approach when it comes to the relevant features of the base topological ring. The key definitions are the following: for a local topological ring $R$ such that $R^\times \subset R$ is open and is a topological group when endowed with the subspace topology, 
\begin{itemize}
\item
$R$ is \emph{\'{e}tale-open} if every \'{e}tale morphism of finite type $R$-schemes is open on $R$-points,

\item
$R$ is \emph{proper-closed} if every proper morphism of finite type $R$-schemes is closed on $R$-points,

\item
$R$ is \emph{finite-closed} if every finite morphism of finite type $R$-schemes is closed on $R$-points.
\end{itemize}
Every local field is \'{e}tale-open and proper-closed (and hence also finite-closed). For further examples as well as a more thorough discussion of these notions, see \S\ref{et-op} and \S\ref{p-c-f-c}.

The axiomatic approach has the added benefit of producing interesting topologies on cohomology beyond the cases when the base is a local field of positive characteristic, or its algebraic extension, or its ring of integers. Other interesting examples in positive characteristic include fields endowed with a Henselian valuation of arbitrary rank, as well as their rings of integers. However, a reader who is only interested in local fields could simplify the discussion by restricting to local field bases throughout.
\epp

\brems
\remi
Philippe Gille and Laurent Moret-Bailly have investigated a similar approach to topologizing cohomology groups and obtained the comparison \Cref{main-Cech} for ground fields $K$ as in \S\ref{p-c-f-c}~\ref{ff-adm-val}. In \S\ref{Hi-top}, we adopt their definition of the topology on $H^n$ for $n \ge 2$ (initially we used the discrete topology for such $n$). This topology is discrete in many cases, see \Cref{sm-discrete}.

\remi \lab{BT-rem}
When all the groups involved are affine, properties of the \v{C}ech topology have also been investigated in \cite{BT14}.\footnote{Op.~cit.~uses a definition of the \v{C}ech topology that a priori differs from the one used here and in \cite{Mil06}. \Cref{BT-counter} shows that the analogue of \cite{BT14}*{Thm.~on p.~562} becomes false for the \v{C}ech topology used here.} Although we have not been able to understand the arguments of op.~cit.~completely, the analogues of \cite{BT14}*{3.1, 4.2.1, 5.1, 5.1.2, and 5.1.3} for the classifying stack topology are special cases of the results of \S\S\ref{topo-coho}--\ref{cont-les}, see Remark \ref{BT-approach} for details. These analogues combine with \Cref{main-Cech,Hi-disc} to give further analogues for the \v{C}ech topology defined in \S\ref{Cech-topo}.
\erems

\bpp[The contents of the paper] 
In \S\ref{app-b}, which may be consulted as needed, we investigate Moret-Bailly's method of topologizing groupoids of rational points of algebraic stacks and also discuss the axiomatic notions mentioned in \S\ref{et-op-ann}. In \S\ref{topo-coho} we study properties of the topologies defined using the classifying stack approach, and in \S\ref{cont-les} we study topological properties of maps arising from a short exact sequence. Although \S\S\ref{topo-coho}--\ref{cont-les} are short, they rely on \S\ref{app-b} and \Cref{BG-genl,new-app-b}, which provide the underlying geometric arguments. \Cref{BG-genl}, included mostly for  convenience, gathers known facts that concern classifying stacks and their interplay with short exact sequences. \Cref{new-app-b} removes diagonal quasi-compactness assumptions from several results in the algebraic space and stack literature. These improvements eliminate a number of quasi-compactness assumptions in \S\S\ref{topo-coho}--\ref{cont-les}; other applications are discussed in \Cref{new-app-b}. The agreement of the classifying stack and the \v{C}ech topologies is the subject of \S\S\ref{comparison}--\ref{comparison-2}: \S\ref{comparison} concerns $H^1$, and \S\ref{comparison-2} deals with $H^n$ for $n \ge 2$.
\epp

\bpp[Notation and conventions] \lab{conv}
For a field $F$, its algebraic closure $\ov{F}$ is chosen implicitly. For a scheme $S$, its fppf site $S_\fppf$ is the category of $S$-schemes with the families of flat, locally of finite presentation, and jointly surjective morphisms as coverings. The small \'{e}tale site is denoted by $S_\et$. The fppf topology is our default choice when taking quotients or considering torsors (synonymously, right torsors). Similarly for cohomology: $H^n$ abbreviates $H^n_\fppf$. Topology on $H^n$ in the absence of the `\v{C}ech' subscript is always that of \S\S\ref{H1-top}--\ref{Hi-top}. For us, compactness does not entail Hausdorffness, and `locally compact' means that every point has a compact neighborhood.

We follow the terminology of \cite{SP} when dealing with algebraic spaces and stacks, albeit for brevity we call a morphism representable if it is representable by algebraic spaces. To emphasize the parallel with schemes, we denote fiber categories by $\sX(S)$ instead of $\sX_S$ and write $x \in \sX(S)$ instead of $x \in \Ob \sX(S)$. We write $\Delta_{\sX/S}$ for the diagonal of an $S$-algebraic stack $\sX$. 
As in \cite{SP}, $\Delta_{\sX/S}$ is neither separated nor quasi-compact at the outset, so algebraic spaces need not be quasi-separated. This is useful: e.g.,~if we were to insist on quasi-compact diagonals, then our results on the topology on $H^1(k, G)$ could not accommodate locally of finite type but not quasi-compact $G$. We ignore set-theoretic difficulties (inherently present in any discussion concerning algebraic stacks) that can be resolved by the use of universes or by the approach used in \cite{SP}.
\epp

\subsection*{Acknowledgements}
I thank Laurent Moret-Bailly for suggesting a number of improvements to our initial results. I thank Bhargav Bhatt, Brian Conrad, Aise Johan de Jong, Bjorn Poonen, Yunqing Tang, and Bertrand To\"{e}n for helpful conversations or correspondence regarding the material of this paper. I thank the referees for helpful comments and suggestions.


\section{Topologizing $R$-points of algebraic stacks} \lab{app-b}

As we recall in \S\ref{sch-case}, the topology on a Hausdorff topological field $K$ gives rise to a natural topology on $X(K)$ for every locally of finite type $K$-scheme $X$. In \cite{MB01}*{\S2}, Moret-Bailly exhibited an elegant way to extend the definition of this topology to the case when $X$ is an algebraic stack. This extension is of major importance for us through the case of a classifying stack $\bbB G$. Therefore, in \S\ref{stack-case} we review the definition of loc.~cit.~(in a slightly more general setting when the base topological ring $R$ need not be a field) and proceed to detail the properties of the resulting topologies in the remainder of \S\ref{app-b}. These properties, especially \Cref{sm-open} \ref{sm-open-a} and \Cref{prop-cl}, are key for our work in \S\S\ref{topo-coho}--\ref{cont-les}.

\bpp[The topology on $R$] \lab{R-def}
Throughout \S\ref{app-b}, $R$ denotes a local topological ring that satisfies
\begin{enumerate}[label=(\greek*)]
\item \lab{alpha}
The group of units $R^\times$ is open in $R$ (equivalently, the maximal ideal $\fm \subset R$ is closed); 

\item \lab{beta}
The inverse map $R^\times \xra{x \mapsto x\i} R^\times$ is continuous when $R^\times \subset R$ is endowed with the subspace topology.
\eenum
\epp

Examples of $R$ that will be of main interest to us are Hausdorff topological fields and arbitrary valuation rings (endowed with their valuation topology).

\bpp[The scheme case] \lab{sch-case}
For locally of finite type $R$-schemes $X$, we want to topologize $X(R)$ so~that
\begin{enumerate}[label=(\roman*)]
\item \lab{cont}
An $R$-morphism $X \ra X\pr$ induces a continuous $X(R) \ra X\pr(R)$;

\item \lab{af}
For each $n\ge 0$, the identification $\bA^n(R) = R^n$ is a homeomorphism;

\item \lab{c-i}
A closed immersion $X \hra X\pr$ induces an embedding $X(R) \hra X\pr(R)$;

\item \lab{o-i}
An open immersion $X \hra X\pr$ induces an open embedding $X(R) \hra X\pr(R)$.
\eenum
\epp
\begin{claim}
Such a topologization must also satisfy
\begin{enumerate}[label=(\roman*)]\addtocounter{enumi}{4}
\item \lab{f-p}
The identifications $(X\pr \times_X X^{\prime\prime})(R) = X\pr(R) \times_{X(R)} X^{\prime\prime}(R)$ are homeomorphisms.
\eenum
\end{claim}

\bpf
The case when $X = \Spec R$ with $X\pr$ and $X^{\prime \prime}$ affine is settled by \ref{af}--\ref{c-i}, which show that any choice of closed immersions $i\pr \colon X' \hra \bA^{n'}$ and $i^{\prime \prime} \colon X^{\prime \prime} \hra \bA^{n^{\prime \prime}}$ induces topological embeddings 
\[
i\pr(R) \colon X'(R) \hra R^{n'}, \q i^{\prime \prime}(R) \colon X^{\prime \prime}(R) \hra R^{n^{\prime \prime}}, \q \text{and} \q (i\pr \times_R i^{\prime \prime})(R) \colon (X' \times_R X^{\prime \prime})(R) \hra R^{n'+n^{\prime \prime}}.
\]

In the case when $X = \Spec R$, we begin by covering $X'$ and $X^{\prime \prime}$ by open affines: $X' = \bigcup_i U'_i$ and $X^{\prime \prime} = \bigcup_j U_j^{\prime \prime}$. Then we recall that $R$ is local to obtain the bijections 
\[
\textstyle X'(R) = \bigcup_i U'_i(R), \q X^{\prime \prime}(R) = \bigcup_j U_j^{\prime \prime}(R), \q \text{and} \q (X' \times_R X^{\prime \prime})(R) = \bigcup_{i, j} (U'_i \times_R U_j^{\prime \prime})(R),
\]
which are covers by open subsets due to \ref{o-i}. We conclude by using the previous case.

In the general case, both the diagonal $\Delta_{X/R}$ and its base change 
\[
X\pr \times_X X^{\prime\prime} \hra X\pr \times_R X^{\prime\prime}
\]
are immersions. Thus, by \ref{c-i}--\ref{o-i}, $(X\pr \times_X X^{\prime\prime})(R)$ has the subspace topology of $(X\pr \times_R X^{\prime\prime})(R)$. However, by the $X = \Spec R$ case, this subspace topology is the topology of $X\pr(R) \times_{X(R)} X^{\prime\prime}(R)$.
\epf

Uniqueness of the sought topologization is ensured by \ref{af}--\ref{o-i}, whereas \cite{Con12b}*{Prop.~3.1} supplies the existence. The conditions \ref{alpha} and \ref{beta} are crucial for the existence---along with the requirement that $R$ is local, they ensure that the topologies on $R$-points interact well with glueing along opens. 

Loc.~cit.~and \cite{Con12b}*{Prop.~2.1} also give further useful properties of the resulting~topologies:
\begin{enumerate}[label=(\roman*)]\addtocounter{enumi}{5}
\item \lab{c-i-H}
If $R$ is Hausdorff, then a closed immersion $X \hra X\pr$ induces a \emph{closed} embedding $X(R)~\hra~X\pr(R)$;

\item \lab{l-c-H}
If $R$ is locally compact and Hausdorff, then $X(R)$ is locally compact.
\eenum

Let $R\pr$ be another local topological ring satisfying \ref{alpha}--\ref{beta} and $R \xra{h} R\pr$ a continuous homomorphism.
\begin{enumerate}[label=(\roman*)]\addtocounter{enumi}{7}
\item \lab{R-Rpr-cts}
Each map $X(R) \ra X(R\pr)$ is continuous.

\item \lab{R-Rpr-o}
If $h$ is an open embedding, then each $X(R) \ra X(R\pr)$ induced by $h$ is open.

\item \lab{R-Rpr-c}
If $h$ is local and a closed embedding, then each $X(R) \ra X(R\pr)$ induced by $h$ is closed. 
\eenum

\bpf[Proof of \ref{R-Rpr-cts}--\ref{R-Rpr-c}]
By considering affine open covers as in the proof of \ref{f-p}, all the claims reduce to the case of an affine $X$. The affine case is in turn addressed in \cite{Con12b}*{Ex.~2.2}.
\epf

\begin{subrem} \lab{R-Rpr-c-rem}
Although loc.~cit.~shows that for affine $X$ the condition that $h$ was local in \ref{R-Rpr-c} is not needed, it cannot be dropped in general. Namely, let $k$ be a nonarchimedean local field, $\fo$ its ring of integers, and $h$ the closed embedding $\fo \hra k$. Build $X$ by glueing the N\'{e}ron lft model of $(\bG_m)_k$ over $\fo$ with $\bA^1_k$ along $(\bG_m)_k \subset \bA^1_k$. Then $X(\fo) \ra X(k)$ is the inclusion $k^\times \hra k$, which is not closed.
\end{subrem}

\bpp[The algebraic space case]
Locally of finite type $R$-algebraic spaces $X$ are in particular algebraic stacks, so we use the procedure described in \S\ref{stack-case} to topologize $X(R)$.
\epp

\bpp[The stack case] \lab{stack-case}
For a locally of finite type $R$-algebraic stack $\sX$, we declare a subcategory $U \subset \sX(R)$ to be \emph{open} if it is full, stable under isomorphisms in $\sX(R)$, and for every $R$-morphism $f\colon X \ra \sX$ with $X$ a locally of finite type $R$-scheme,  $f(R)\i(U)$ is open in $X(R)$ (the latter is topologized as in \S\ref{sch-case}). In the last condition, due to \ref{o-i}, one can also restrict to affine $X$. 

The above definition topologizes $\sX(R)$ in the sense that arbitrary unions and finite intersections of opens are open. In particular, the set of isomorphism classes of objects of $\sX(R)$ becomes a bona fide topological space, which, thanks to \ref{cont}, is none other than that of \S\ref{sch-case} if $\sX$ is representable by a scheme. The topology on this set of isomorphism classes uniquely recovers the collection of open subcategories, and hence also the ``topology'', of $\sX(R)$. When discussing the latter we freely use evident analogues of familiar notions (e.g.,~closedness or continuity), which always correspond to bona fide notions after taking isomorphism classes of objects. 

Some of the properties \ref{cont}--\ref{R-Rpr-c} have analogues for algebraic stacks, see \Cref{basic} and \Cref{sm-open} \ref{st-dir-pr}, \ref{sm-open-b}, and \ref{sm-open-c}.
\epp

\bpp[The presheaf case] \lab{presheaf-case}
For a set-valued presheaf $\sF$ on the category of locally of finite type $R$-schemes, we declare $U \subset \sF(R)$ to be \emph{open} if for every presheaf morphism $f\colon X \ra \sF$ with $X$ a locally of finite type $R$-scheme, $f(R)\i(U)$ is open in $X(R)$. Similarly to \S\ref{stack-case}, this topologizes $\sF(R)$. A presheaf morphism $\sF \ra \sF\pr$ induces a continuous $\sF(R) \ra \sF\pr(R)$. If $\sF$ is the presheaf of isomorphism classes of objects of fiber categories of a locally of finite type $R$-algebraic stack $\sX$, then the resulting topology on $\sF(R)$ agrees with that of \S\ref{stack-case}, and hence also with that of \S\ref{sch-case} if $\sX$ is a scheme (i.e., if $\sF$ is representable). 

Our main case of interest is $\sF = H^n(-, G)$ for an $n \in \bZ_{\ge 0}$ and a commutative flat locally of finite presentation $R$-group algebraic space $G$; in this case, see \Cref{sm-discrete} \ref{disc-2}--\ref{disc-4} and \Cref{top-gp} for conditions that ensure that $\sF(R)$ is a topological group.
\epp

In the study of openness or closedness of induced maps on $R$-points, the following lemma often provides a robust reduction to the algebraic space case.

\blem \lab{op-cl}
For a representable morphism $f\colon \sX \ra \sY$ of locally of finite type $R$-algebraic stacks and an $R$-morphism $s\colon Y \ra \sY$ with $Y$ a finite type affine $R$-scheme, consider the Cartesian square
\be\ba \lab{op-cl-eq}
\xymatrix{
X \ar[d]_-{s\pr}\ar[r]^-{f\pr} & Y \ar[d]_{s} \\
\sX \ar[r]^-{f} & \sY.
}
\ea\ee
If $f\pr(R)$ is open (resp.,~closed) for every $s$ as above, then $f(R)$ is open (resp.,~closed).
\elem

\bpf
For an open (resp.,~closed) $U \subset \sX(R)$, its full essential image $f(R)(U) \subset \sY(R)$ is open (resp.,~closed) if so is $s(R)\i(f(R)(U))$ for every $s$. It remains to observe the equality 
\[
s(R)\i(f(R)(U)) = f\pr(R)(s\pr(R)\i(U))
\]
that is provided by the Cartesian property of \eqref{op-cl-eq}.
\epf

\bcor \lab{basic}
The topologies of \S\ref{stack-case} satisfy the analogues of \ref{cont}, \ref{o-i}, \ref{c-i-H}, and \ref{R-Rpr-cts}.
\ecor

\bpf
\Cref{op-cl} gives \ref{o-i} and \ref{c-i-H}; \ref{R-Rpr-cts} follows from the scheme case; \ref{cont} is clear.
\epf

\bpp[\'{E}tale-openness] \lab{et-op}
Let $R$ be a local topological ring that satisfies \ref{alpha}--\ref{beta}. We call $R$ \emph{\'{e}tale-open} if for every \'{e}tale morphism $f\colon X \ra Y$ of locally of finite type $R$-schemes, $f(R)$ is open. Due to \ref{o-i}, the further requirement that $X$ and $Y$ are affine results in the same class of \'{e}tale-open $R$. 

Before proceeding we give some examples of \'{e}tale-open $R$; all $R$ in the examples are Henselian.
\epp
\begin{enumerate}[label={(\arabic*)}] 
\item \lab{R-and-C}
$\bR$ and $\bC$.

\item \lab{ff-Hens-val}
The fraction field $K$ of a Henselian valuation ring $A$.
\eenum
\bpf
If $A = K$, then $K$ and $Y(K)$ are discrete. If $A \neq K$, then, according to \cite{GGMB14}*{3.1.4}, for every \'{e}tale morphism $f\colon X \ra Y$ between finite type $K$-schemes, the induced morphism $f(K)\colon X(K) \ra Y(K)$ is a local homeomorphism; in particular, $f(K)$ is open.
\epf
\begin{enumerate}[label={(\arabic*)}] \addtocounter{enumi}{2}
\item \lab{Hens-val}
A Henselian valuation ring $A$.
\eenum
\bpf
We use \ref{o-i} to assume that $X$ and $Y$ are affine, and we let $K$ denote the fraction field of $A$. By \ref{R-Rpr-cts}--\ref{R-Rpr-o}, both $X(A) \hra X(K)$ and $Y(A) \hra Y(K)$ are inclusions of open subsets, so the openness of $f(A) \colon X(A) \ra Y(A)$ results from that of $f(K) \colon X(K) \ra Y(K)$ supplied by \ref{ff-Hens-val}.
\epf
\begin{enumerate}[label={(\arabic*)}] \addtocounter{enumi}{3}
\item \lab{nonarch}
A nonarchimedean local field and its ring of integers.
\eenum
\bpf
These are special cases of \ref{ff-Hens-val} and \ref{Hens-val}.
\epf

Due to the \'{e}tale-local nature of algebraic spaces and the smooth-local nature of algebraic stacks, the \'{e}tale-openness condition is key for the topologies on $R$-points of algebraic stacks to be well-behaved. This condition is often accompanied by the requirement that $R$ is Henselian, which through \Cref{lift-a-lot} ensures that $R$-points lift to some smooth scheme neighborhood (such liftings facilitate reductions to scheme cases). We record some of the resulting desirable properties in the following proposition.

\bprop \lab{sm-open}
Let $h \colon R \ra R\pr$ be a continuous homomorphism between Henselian local topological rings that satisfy \ref{alpha}--\ref{beta} and are \'{e}tale-open, and let $\sX$, $\sY$ be locally of finite type $R$-algebraic stacks such that $\Delta_{\sX/R}$ has a separated $R$-fiber over the closed point of $\Spec R$. 
\benum
\item \lab{sm-open-a}
If $f\colon \sX \ra \sY$ is a smooth morphism, then $f(R)$ is open. 

\item \lab{st-dir-pr}
If $\Delta_{\sY/R}$ also has a separated $R$-fiber over the closed point of $\Spec R$, then the identification 
\[
(\sX \times_R \sY)(R) \isomto \sX(R) \times \sY(R)
\]
is a homeomorphism (compare with~\ref{f-p}). 

\item \lab{sm-open-b}
If $R$ is locally compact and Hausdorff, then $\sX(R)$ is locally compact (compare with \ref{l-c-H}).

\item \lab{compact}
If $R$ is compact Hausdorff and the number of isomorphism classes of objects of $\sX(R/\fm)$ is finite, where $\fm \subset R$ is the maximal ideal, then $\sX(R)$ is compact.

\item \lab{sm-open-c}
If $h$ is an open embedding, then the map 
\[
\sX(R) \ra \sX(R\pr)
\]
induced by $h$ is open (compare with \ref{R-Rpr-o}).
\eenum
\eprop

\bpf
We use \Cref{lift-a-lot} to choose a smooth $\sX$-scheme $X$ with $X(R) \ra \sX(R)$ essentially surjective; for \ref{compact}, we choose $X$ to be, in addition, affine.

\benum
\item 
We use the continuity and the essential surjectivity of $X(R) \ra \sX(R)$ to replace $\sX$ by $X$, and then we use \Cref{op-cl} to assume that $\sY$ is a scheme and $\sX$ is an algebraic space. We then once more replace $\sX$ by $X$ to reduce to the case when both $\sX$ and $\sY$ are schemes.

Due to \ref{o-i}, we may work locally on $\sX$. Therefore, due to \cite{BLR90}*{\S2.2/11}, we may assume that $f$ factors as $\sX \xra{g} \bA^n_\sY \xra{\pi} \sY$ where $g$ is \'{e}tale and $\pi$ is the canonical projection. It remains to observe that $g(R)$ is open by the \'{e}tale-openness of $R$ and that $\pi(R)$ is a projection onto a direct factor by \ref{f-p}, and hence is also open.

\item
We use \Cref{lift-a-lot} to choose a smooth $\sY$-scheme $Y$ with $Y(R) \ra \sY(R)$ essentially surjective. By the aspect \ref{cont} of \Cref{basic}, the identification in question is continuous. By \ref{sm-open-a}, it inherits openness from the homeomorphism 
\[
\qq (X \times_R Y)(R) \isomto X(R) \times Y(R).
\]

\item
The claim follows by combining \ref{sm-open-a} applied to $X \ra \sX$ with \ref{l-c-H} applied to $X$.

\item
Since $X$ is affine, $X(R)$ is compact by \ref{af} and~\ref{c-i-H}; thus, so is its continuous image $\sX(R)$.

\item
By \ref{R-Rpr-o} applied to $X$, the map $X(R) \ra X(R')$ induced by $h$ is open. By \ref{sm-open-a}, the maps 
\[
\qq X(R) \ra \sX(R) \qq \text{and} \qq X(R') \ra \sX(R')
\]
are also open; moreover, $X(R) \ra \sX(R)$ is essentially surjective by construction. Therefore, by combining these observations, the map $\sX(R) \ra \sX(R\pr)$ induced by $h$ is open, too. \qedhere

\eenum
\epf

\bpp[Restriction of scalars] \lab{rest-scal}
Let $R$ and $R\pr$ be local topological rings that satisfy \ref{alpha}--\ref{beta} and let $h\colon R \ra R\pr$ be a continuous homomorphism that makes $R\pr$ a finite free $R$-module (both algebraically and topologically). Let $\sX\pr$ be a locally of finite type $R\pr$-algebraic stack and let 
\[
\sX \ce \Res_{R\pr/R} \sX\pr
\]
be its restriction of scalars defined as the category fibered in groupoids over $\Spec R$ with $S$-points $\sX\pr(S_{R\pr})$, functorially in the $R$-scheme $S$. Assume that $\sX$ is a locally of finite type $R$-algebraic stack. By \cite{Ols06}*{1.5}, this assumption is met if $\sX'$ is locally of finite presentation over $R'$ and $\Delta_{\sX'/R'}$ is finite; by \cite{SP}*{\href{http://stacks.math.columbia.edu/tag/05YF}{05YF} and \href{http://stacks.math.columbia.edu/tag/04AK}{04AK}} (compare with \cite{BLR90}*{7.6/4--5}), it is also met if $\sX'$ is an algebraic space that is locally of finite presentation over $R'$. 

We wish to describe some situations in which the identification $\sX'(R') \cong \sX(R)$ is a homeomorphism.
\epp

\bprop \lab{weil-weil}
In \S\ref{rest-scal}, the identification 
\[
\iota\colon \sX\pr(R\pr) \cong \sX(R)
\]
is a homeomorphism
\benum
\item \lab{weil-weil-a}
If $\sX\pr$ is a scheme, or

\item \lab{weil-weil-b}
If $R$ and $R\pr$ are Henselian \'{e}tale-open and $\Delta_{\sX\pr/R\pr}$ has a separated $R\pr$-fiber over the closed point of $\Spec R\pr$. 
\eenum
\eprop

\bpf Let $x\pr \in \sX\pr(R\pr)$ be arbitrary.
\benum
\item 
If $\sX\pr$ is affine, then let $\sX\pr \hra \bA^n_{R\pr}$ be a closed immersion. By \cite{BLR90}*{7.6/2 (ii)}, the restriction of scalars $\sX \hra \bA^{dn}_{R}$ with $d$ determined by $R\pr \cong R^{\oplus d}$ is also a closed immersion. Therefore, since the isomorphism $R\pr \cong R^{\oplus d}$ agrees with the topologies, \ref{af} and \ref{c-i} settle the case of affine $\sX\pr$.  

If $\sX\pr$ is arbitrary, then let $U' \subset \sX'$ be an affine open through which $x'$ factors. By \cite{BLR90}*{7.6/2 (i)}, the restriction of scalars $U \subset \sX$ is an open immersion. Therefore, by the aspect \ref{o-i} of \Cref{basic}, 
\[
\qq U(R) \subset \sX(R) \qq \text{and} \qq U'(R') \subset \sX'(R')
\]
are inclusions of open subsets, so $\iota$ is a local homeomorphism at $x'$ by the settled affine case. Since $\iota$ is bijective and $x'$ was arbitrary, $\iota$ must be a homeomorphism.

\item
By \Cref{LMB-revamp}, $x\pr$ lifts to an $\wt{x\pr} \in X\pr(R\pr)$ for some affine scheme $X\pr$ equipped with a smooth morphism $f\pr\colon X\pr \ra \sX\pr$. By the infinitesimal lifting criterion \cite{LMB00}*{4.15},\footnote{The proof of loc.~cit.~makes no use of the general conventions of \cite{LMB00} that the base schemes are quasi-separated and that the diagonals of algebraic stacks are separated and quasi-compact. } 
\[
\qq \Res_{R\pr/R} f\pr\colon X \ra \sX
\]
is smooth. Thus, \ref{weil-weil-a} and \Cref{sm-open} \ref{sm-open-a} show that $\iota$ is a local homeomorphism at $x'$. Since $\iota$ is also bijective, it must be a homeomorphism. \qedhere
\eenum
\epf

In addition to \'{e}tale-openness, another pleasant situation is when finite or, more generally, proper morphisms induce closed maps on $R$-points. We axiomatize this in the following notions.

\bpp[Proper-closedness and finite-closedness] \lab{p-c-f-c}
Let $R$ be a local topological ring that satisfies \ref{alpha}--\ref{beta}. We call $R$ \emph{finite-closed} (resp.,~\emph{proper-closed}) if for every finite (resp.,~proper) morphism $f\colon X \ra Y$ of locally of finite type $R$-schemes, $f(R)$ is closed. Thanks to \ref{o-i}, these concepts are unaltered if one further requires $Y$ to be affine. A proper-closed $R$ is also finite-closed. An $R$ that is finite-closed (or proper-closed) is also Hausdorff: take $f$ to be the inclusion of the origin of $\bA^1_R$.
\epp

\begin{lem} \lab{Pn-cpct}
If $k$ is a local field, then 
\benum
\item \lab{Pn-cpct-a}
The field $k$ is proper-closed (and hence also finite-closed);

\item \lab{Pn-cpct-b}
The topological space $X(k)$ is compact for every proper $k$-algebraic space $X$;

\item \lab{Pn-cpct-d}
The map $f(k)$ is closed for every proper representable morphism $f\colon \sX \ra \sY$ of locally of finite type $k$-algebraic stacks. 
\eenum
If, in addition, $k$ is nonarchimedean and $\fo$ is its ring of integers, then \ref{Pn-cpct-a} and \ref{Pn-cpct-b} also hold with $\fo$ in place of $k$ (for \ref{Pn-cpct-d}, see \Cref{turbo-c}). \end{lem}
\bpf
The key input is Chow's lemma due to Knutson \cite{Knu71}*{IV.3.1}; we use its following form:

\begin{lem}[\cite{SP}*{\href{http://stacks.math.columbia.edu/tag/088U}{088U}}] \lab{Chow}
For a proper map $f\colon X \ra Y$ of separated $k$-algebraic spaces of finite type, there is a proper map $g\colon X\pr \ra X$ of $k$-algebraic spaces such that $g\i(U) \cong U$ for a dense open subspace $U \subset X$ and $f \circ g$ factors through a closed immersion into $\bP^n_Y$ for some $n \ge 0$. 
\end{lem}

\Cref{Chow} suggests the following strategy, which we borrow from \cite{Con12b}*{proof of Prop.~4.4}: reduce the claim in question to the cases when $X$ is replaced by a relative projective space or the reduced closed subspace $Z \subsetneq X$ complementary to $U$. Since the underlying topological space of $X$ is Noetherian, iteration will shrink $Z$ to the empty space, leaving only the first case to consider.
\benum
\item 
Assume that $Y$ in \Cref{Chow} is an affine scheme. Since 
\be \lab{pts-union}
X(k) = g(k)(X\pr(k)) \union Z(k),
\ee
the strategy above and \ref{f-p}--\ref{c-i-H} reduce closedness of $f(k)$ to that of the second projection 
\[
\bP^n(k) \times Y(k) \ra Y(k) \qq \text{for $n \ge 0$,}
\]
which follows from the compactness of~$\bP^n(k)$.

\item
Combine the strategy above, \eqref{pts-union}, and the compactness of $X\pr(k)$ inherited from $\bP^n(k)$.

\item
Combine \Cref{op-cl} and the proof of \ref{Pn-cpct-a}. 
\eenum

We now turn to the variants for $\fo$ in place of $k$. All of them reduce to their counterparts for $k$.
\benum
\item \lab{o-pf}
Fix a proper morphism $f\colon X \ra Y$ of finite type $\fo$-schemes with $Y$ affine, so $Y(\fo) \subset Y(k)$ is a closed embedding by \ref{R-Rpr-c} and \Cref{R-Rpr-c-rem}. By the valuative criterion, 
\[
\qq f(k)\i(Y(\fo)) = X(\fo)
\]
in $X(k)$, so $f(k)|_{X(\fo)} = f(\fo)$ by \ref{R-Rpr-o}.

\item \lab{o-pf-2}
By the valuative criterion of properness for algebraic spaces \cite{SP}*{\href{http://stacks.math.columbia.edu/tag/0A40}{0A40}}, \Cref{basic}~\ref{R-Rpr-cts}, and \Cref{sm-open} \ref{sm-open-c} (with \S\ref{et-op}~\ref{Hens-val}), $X(\fo) \ra X(k)$ is a homeomorphism. If $X$ is a scheme, \ref{R-Rpr-o} applies directly.  \qedhere
\eenum
\epf

Before proceeding we give some examples of proper-closed and finite-closed $R$:

\begin{enumerate}[label={(\arabic*)}]
\item \lab{loc-prop-cl} 
A local field $k$ is proper-closed. So is its ring of integers $\fo$ if $k$ is nonarchimedean.
\eenum

\bpf
This is a special case of \Cref{Pn-cpct}.
\epf

\begin{enumerate}[label={(\arabic*)}] \addtocounter{enumi}{1}
\item \lab{ff-adm-val}
Let $A$ be a Henselian valuation ring, $K \ce \Frac(A)$, and $\wh{K}$ the completion of $K$. If $\wh{K}/K$ is a separable field extension, then both $A$ and $K$ are finite-closed.
\eenum

\bpf
If $A = K$, then $K$ and $Y(K)$ are discrete. If $A \neq K$, then $K$ is finite-closed by \cite{GGMB14}*{4.2.6}, and $A$ is then also finite-closed by the method of proof of \Cref{Pn-cpct} \ref{Pn-cpct-a} for $\fo$. 
\epf



\brem
The \'{e}tale-open and finite-closed conditions have already appeared in the literature: they are properties (H) and (F) in \cite{MB12}*{\S2.1} in the case of a Hausdorff topological~field.
\erem

In \Cref{prop-cl,gamma,turbo-c} we record several situations in which morphisms of algebraic stacks (e.g.,~of schemes) induce closed maps on $R$-points. These results will be useful in \S\ref{cont-les}, see the proof of \Cref{cl-supreme}. \Cref{prop-cl} allows quite general $R$ at the expense of stringent conditions on the morphism, whereas \Cref{gamma,turbo-c} allow more general morphisms but restrict $R$.

\bprop \lab{prop-cl}
Let $R$ be a local topological ring that satisfies \ref{alpha}--\ref{beta}, and let 
\[
f\colon \sX \ra \sY
\]
be a morphism of locally of finite type $R$-algebraic stacks.
\benum
\item \lab{prop-cl-b}
If $R$ is finite-closed and $f$ is finite (and hence representable by schemes), then $f(R)$ is a closed map.

\item \lab{prop-cl-a}
If $R$ is proper-closed and $f$ is proper and representable by schemes, then $f(R)$ is a closed map.
\eenum
\eprop

\bpf
The claims follow from \Cref{op-cl} and from the definitions.
\epf

\bprop \lab{gamma}
Let $R$ be a local topological ring that is Hausdorff, satisfies \ref{alpha}--\ref{beta}, and such that $R^\times$ is closed in $R$ (e.g., $R$ could be a valuation ring). Let 
\[
f\colon \sX \ra \sY
\]
be a representable by schemes morphism of locally of finite type $R$-algebraic stacks.
\benum
\item \lab{gamma-a}
If $f$ is a quasi-compact immersion, then $f(R)$ is a closed embedding.

\item \lab{gamma-b}
If $R$ is finite-closed and $f$ is quasi-finite, then $f(R)$ is a closed map.

\item \lab{gamma-c}
If $R$ is proper-closed and $f$ is of finite type, then $f(R)$ is a closed map. 
\eenum
\eprop

\bpf
In \ref{gamma-a}, we need to show that the continuous monomorphism $f(R)$ is closed. In all cases \Cref{op-cl} permits us to assume that $\sY$ is an affine scheme $\Spec B$ (so $f$ is a scheme morphism).
\benum
\item
By \S\ref{sch-case} \ref{c-i-H}, we may assume that $f$ is a quasi-compact open immersion. By quasi-compactness, we may further assume that $f$ is the inclusion of a principal open affine $\Spec B[\f{1}{b}]$. We view $b$ as a morphism $b\colon \sY \ra \bA^1_R$. Then the closedness of $R^\times$ implies that of 
\[
\qq f(R)(\sX(R)) = b(R)\i(R^\times),
\]
so the claim results from \S\ref{sch-case} \ref{o-i}.
\eenum

Thanks to the quasi-compactness of $f$, in the proof of \ref{gamma-b} and \ref{gamma-c} we may replace $\sX$ by an open affine, and hence further assume that $f$ is separated.

\benum \addtocounter{enumi}{1}

\item
Combine \cite{EGAIV4}*{18.12.13} (i.e.,~Zariski's main theorem) and \ref{gamma-a}.

\item
Combine \cite{Con07b}*{Thm.~4.1} (i.e.,~Nagata's embedding theorem) and \ref{gamma-a}. \qedhere
\eenum
\epf

\bprop \lab{turbo-c}
Let $\fo$ be the ring of integers of a nonarchimedean local field $k$. For every representable, separated, finite type morphism 
\[
f\colon \sX \ra \sY
\] 
of locally of finite type $\fo$-algebraic stacks, the map $f(\fo)$ is closed.
\eprop

\bpf
\Cref{op-cl} reduces to the case when $\sY$ is an affine scheme. Nagata's embedding theorem for algebraic spaces \cite{CLO12}*{Thm.~1.2.1} and \Cref{gamma} \ref{gamma-a} reduce further to proper $f$. By \S\ref{sch-case}~\ref{R-Rpr-cts} and the first part of the first sentence of \Cref{R-Rpr-c-rem}, $\sY(\fo) \subset \sY(k)$ is a closed embedding. By the valuative criterion of properness for algebraic spaces \cite{SP}*{\href{http://stacks.math.columbia.edu/tag/0A40}{0A40}}, $f(k)\i(\sY(\fo)) = \sX(\fo)$ in $\sX(k)$, so the embedding $\sX(\fo) \subset \sX(k)$ that results from \Cref{sm-open} \ref{sm-open-c} has a closed image, and hence is a closed embedding. It remains to note that $f(k)$ is closed by \Cref{Pn-cpct}~\ref{Pn-cpct-d}.
\epf

We close \S\ref{app-b} by surveying conditions which ensure that $\sX(R)$ is Hausdorff.

\bprop\lab{H-supreme}
Let $R$ be a local topological ring that is Hausdorff and satisfies \ref{alpha}--\ref{beta}, and let $\sX$ be a locally of finite type $R$-algebraic stack such that $\Delta_{\sX/R}$ has a separated $R$-fiber over the closed point of $\Spec R$. If 
\begin{enumerate}[label={(\arabic*)}]
\item \lab{H-1}
$\sX$ is a separated $R$-scheme, or

\item \lab{H-2}
$R^\times$ is closed in $R$ (e.g.,~$R$ is valuation ring), and $\sX$ is a quasi-separated $R$-scheme, 
\eenum
then $\sX(R)$ is Hausdorff. If $R$ is Henselian and \'{e}tale-open and
\begin{enumerate}[label={(\arabic*)}] \addtocounter{enumi}{2}

\item \lab{H-3}
$\sX$ is a separated $R$-algebraic space, or

\item \lab{H-4}
$R^\times$ is closed in $R$, and $\sX$ is a quasi-separated $R$-algebraic space with $\Delta_{\sX/R}$ an immersion,~or

\item \lab{H-5}
$R$ is finite-closed and $\Delta_{\sX/R}$ is finite, or

\item \lab{H-6}
$R$ is finite-closed, $R^\times$ is closed in $R$, and $\Delta_{\sX/R}$ is quasi-finite and separated, or

\item \lab{H-9}
$R$ is a local field and $\sX$ is $R$-separated (i.e.,~$\Delta_{\sX/R}$ is proper), or

\item \lab{H-10}
$R$ is the ring of integers of a nonarchimedean local field, and $\Delta_{\sX/R}$ is quasi-compact and separated, or 

\item \lab{H-7}
$R$ is proper-closed, and $\Delta_{\sX/R}$ is proper and representable by schemes, or

\item \lab{H-8}
$R$ is proper-closed, $R^\times$ is closed in $R$, and $\Delta_{\sX/R}$ is quasi-compact and representable by~schemes,
\eenum
then $\sX(R)$ is Hausdorff. 
\eprop

\bpf
The identification $(\sX \times_R \sX)(R) \isomto \sX(R) \times \sX(R)$ is a homeomorphism: 

\q\ In the cases \ref{H-1}--\ref{H-2}, this follows from \S\ref{sch-case} \ref{f-p}; 

\q\ In the cases \ref{H-3}--\ref{H-8}, this follows from \Cref{sm-open} \ref{st-dir-pr}. 

It remains to note that the diagonal map $\Delta_{\sX/R}(R)$ is closed:

\q\ In the case \ref{H-1}, this follows from \S\ref{sch-case} \ref{c-i-H}; 

\q\ In the cases \ref{H-2} and \ref{H-4}, this follows from \Cref{gamma} \ref{gamma-a}; 

\q\ In the case \ref{H-3}, this follows from the aspect \ref{c-i-H} of \Cref{basic}; 

\q\ In the case \ref{H-5}, this follows from \Cref{prop-cl} \ref{prop-cl-b}; 

\q\ In the case \ref{H-6}, this follows from \Cref{gamma} \ref{gamma-b} (representability of $\Delta_{\sX/R}$ by schemes is guaranteed by \cite{LMB00}*{A.2}); 

\q\ In the case \ref{H-9}, this follows from \Cref{Pn-cpct} \ref{Pn-cpct-d}; 

\q\ In the case \ref{H-10}, this follows from \Cref{turbo-c} (quasi-compactness ensures that $\Delta_{\sX/R}$ is of finite type because it is always locally of finite type by \cite{SP}*{\href{http://stacks.math.columbia.edu/tag/04XS}{04XS}} or by \cite{LMB00}*{4.2 and its proof}); 

\q\ In the case \ref{H-7}, this follows from \Cref{prop-cl} \ref{prop-cl-a}; 

\q\ In the case \ref{H-8}, this follows from \Cref{gamma} \ref{gamma-c} (supplemented by \cite{SP}*{\href{http://stacks.math.columbia.edu/tag/04XS}{04XS}} or by \cite{LMB00}*{4.2} again). \qedhere
\epf


\section{Topologies on cohomology sets via classifying stacks} \lab{topo-coho}

The setup introduced in \S\ref{H1-top} is valid throughout \S\ref{topo-coho}. Our main goal in \S\S\ref{topo-coho}--\ref{cont-les} is to detail the properties of the topologies introduced in \S\S\ref{H1-top}--\ref{Hi-top}.

\bpp[Topology on $H^0(R, G)$ and $H^1(R, G)$] \lab{H1-top}
Let $R$ be a local topological ring that satisfies \ref{alpha} and \ref{beta} of \S\ref{R-def}, i.e., $R^\times \subset R$ is open and a topological group when endowed with the subspace topology. Let $G$ be a flat and locally of finite presentation $R$-group algebraic space, and let $\bbB G$ be its classifying stack, which is algebraic and $R$-smooth by \S\ref{BG-alg} and \Cref{BG-sm}. We use \S\ref{stack-case} to topologize 
\[
H^0(R, G) = G(R) \qq \text{and} \qq H^1(R, G)
\]
by viewing the latter as the set of isomorphism classes of objects of $\bbB G(R)$. 

Although in most examples $G$ is a scheme, allowing algebraic spaces leads to a more robust theory and more flexible proofs: for instance, group algebraic spaces behave better than group schemes with respect to representability of quotients or of inner forms.
\epp

\bpp[Topology on $H^n(R, G)$ for $n \ge 2$] \lab{Hi-top}
If in the setup of \S\ref{H1-top}, $G$ is in addition commutative, then we use the definitions of \S\ref{presheaf-case} with $\sF = H^n(-, G)$ to topologize 
\[
H^n(R, G) \qq \text{for $n \ge 2$.}
\]
(One can also use \S\ref{presheaf-case} in the $n = 0$ and $n = 1$ cases---the discussion in \S\ref{presheaf-case} guarantees agreement with \S\ref{H1-top}.) 
\epp

The following two propositions illustrate the functoriality inherent in the definitions.

\bprop\lab{map-cts}
For a homomorphism 
$G \ra G'$
between flat and locally of finite presentation $R$-group algebraic spaces, the induced maps
\[
H^n(R, G) \ra H^n(R, G')
\]
are continuous for $n \le 1$, and also for $n \ge 2$ if $G$ and $G'$ are commutative.
\eprop

\bpf
There is an underlying presheaf morphism 
\[
H^n(-, G) \ra H^n(-, G'),
\]
so the continuity follows from the discussion of \S\ref{presheaf-case}.
\epf

\bprop \lab{twist}
Let $T$ be a right $G$-torsor and let 
${}_TG \ce \Aut_G(T)$
be the corresponding inner form of $G$. Twisting by $T$ induces a homeomorphism 
\[
H^1(R, G) \cong H^1(R, {}_TG)
\]
that sends the class of $T$ to the neutral class.
\eprop

\bpf
By \cite{Gir71}*{III.2.6.1 (i)}, twisting by $T$ induces an underlying isomorphism 
\[
\bbB G \cong \bbB{}_TG
\]
of algebraic stacks over $R$. This isomorphism is defined by $X \mapsto \Hom_G(T, X)$, and hence sends the class of $T$ to the neutral class.
\epf

Before proceeding to discuss finer topological properties, we record conditions which ensure discreteness of the topologies introduced in \S\S\ref{H1-top}--\ref{Hi-top}.

\bprop \lab{sm-discrete}
If $R$ is Henselian and \'{e}tale-open (in the sense of \S\ref{et-op}), then 
\benum 
\item \lab{disc-1}
$H^1(R, G)$ is discrete if $G$ is smooth;

\item \lab{disc-2}
$H^n(R, G)$ is discrete for $n \ge 1$ if $G$ is commutative and smooth;

\item \lab{disc-3}
$H^n(R, G)$ is discrete for $n \ge 2$ if $G$ is commutative and of finite presentation;

\item \lab{disc-4}
$H^n(R, G)$ is discrete for $n \ge 2$ if $G$ is commutative and $R$ is excellent (and hence Noetherian).
\eenum
\eprop

\bpf 
We fix a locally of finite type $R$-scheme $X$, an $a \in H^n(R, G)$, and an $f \in H^n(X, G)$, which corresponds to a presheaf morphism $f\colon X \ra H^n(-, G)$. We need to show that 
\[
f(R)\i(a) \subset X(R)
\]
is open. For every $x\in f(R)\i(a)$, we will build an \'{e}tale $g_x\colon U_x \ra X$ such that 
\[
x \in g_x(R)(U_x(R)) \qq \text{and} \qq f|_{U_x} = a|_{U_x} \ \text{in $H^n(U_x, G)$.} 
\]
The latter forces 
\[
g_x(R)(U_x(R)) \subset f(R)\i(a),
\]
so the \'{e}tale-openness of $R$ will prove the desired openness of $f(R)\i(a) \subset X(R)$.

For the construction of $U_x$, we view $x$ as an $R$-morphism $\Spec R \xra{x} X$, let $s \in \Spec R$ be the closed point, and let $\cO$ be the Henselization of $X$ at $x(s)$. Since $R$ is Henselian, $x$ factors through $\Spec \cO \ra X$; being an $R$-morphism, $x$ also identifies the residue fields at $s$ and $x(s)$. Therefore, 
\[
f|_{\cO} = a|_{\cO} \ \ \text{in $H^n(\cO, G)$}
\]
by Remark \ref{sm-H1-inj}  for \ref{disc-1}, \Cref{Gro68-revamp} for \ref{disc-2}, \Cref{Toe11-revamp} for \ref{disc-3}, and \cite{Toe11}*{3.4} together with \cite{EGAIV2}*{7.8.3} and \cite{EGAIV4}*{18.7.6} for \ref{disc-4}. It remains to employ limit formalism \cite{SGA4II}*{VII, 5.9}\footnote{We use analogues of loc.~cit.~for algebraic spaces instead of schemes and fppf cohomology of algebraic spaces instead of \'{e}tale cohomology of schemes. The proofs based on topos-theoretic generalities of \cite{SGA4II}*{VI} are analogous,~too.} to find a desired $g_x$ through which $\Spec \cO \ra X$ factors.
\epf

In the remainder of \S\ref{topo-coho} we suppress mentioning $H^n(R, -)$ for $n \ge 2$---in the cases of such $n$ we have nothing to add to the frequent discreteness given by \Cref{sm-discrete}~\ref{disc-2}--\ref{disc-4}.

For $n = 0$ and $n = 1$, the following conditions ensure that $H^n(R, G)$ is a topological group.

\bprop \lab{top-gp} Let $\fm \subset R$ be the maximal ideal.
\benum
\item \lab{top-gp-naive} 
If $G$ is a scheme, then $H^0(R, G)$ is a topological group.

\item \lab{top-gp-H0} 
If $R$ is Henselian and \'{e}tale-open, then $H^0(R, G)$ is a topological group.

\item \lab{top-gp-H1} 
If $R$ is Henselian and \'{e}tale-open, $G$ is commutative, and either $G$ is smooth or $G_{R/\fm}$ is a scheme (by \cite{Art69a}*{4.2}, $G_{R/\fm}$ is a scheme if, for instance, it is quasi-separated), then $H^1(R, G)$ is an abelian topological group.
\eenum
\eprop

\bpf
The inverse map of the group $H^n(R, G)$ is continuous because it is induced by an underlying morphism of algebraic stacks (even of schemes in \ref{top-gp-naive} and of algebraic spaces in \ref{top-gp-H0}), as mentioned in the last paragraph of \S\ref{BG-alg} in the case \ref{top-gp-H1}. Likewise, we deduce that the multiplication map is also continuous once we argue that the bijection
\be\lab{top-gp-bij}
H^n(R, G\times G) \cong H^n(R, G) \times H^n(R, G)
\ee
is a homeomorphism. 

\q\ In \ref{top-gp-naive}, \eqref{top-gp-bij} is a homeomorphism by \S\ref{sch-case} \ref{f-p}. 

\q\ In \ref{top-gp-H0}, \eqref{top-gp-bij} is a homeomorphism by \Cref{sm-open} \ref{st-dir-pr}.

\q\ In \ref{top-gp-H1}, if $G$ is smooth, then both sides of \eqref{top-gp-bij} are discrete by \Cref{sm-discrete} \ref{disc-1}.

\q\ In \ref{top-gp-H1}, if $G_{R/\fm}$ is a scheme, then $G_{R/\fm}$ is separated, so $\Delta_{\bbB G_{R/\fm}}$ is also separated by \Cref{Delta-BG}~\ref{Delta-BG-b}, and hence \eqref{top-gp-bij} is a homeomorphism by \Cref{sm-open} \ref{st-dir-pr}.
\epf

In the following proposition we record conditions which ensure that $H^n(R, G)$ is locally compact. Local compactness (as well as Hausdorffness) is important to know in practice: for instance, Pontryagin duality concerns locally compact Hausdorff abelian topological groups, so to make sense of Tate--Shatz local duality over nonarchimedean local fields of positive characteristic one needs to know that the cohomology groups in question are locally compact and Hausdorff.

\bprop 
Suppose that $R$ is locally compact and Hausdorff. 
\benum
\item
If $G$ is a scheme, then $H^0(R, G)$ is locally compact.

\item
If $R$ is Henselian and \'{e}tale-open, then $H^0(R, G)$ is locally compact.

\item 
If $R$ is Henselian and \'{e}tale-open and $G_{R/\fm}$ is a scheme, then $H^1(R, G)$ is locally compact.
\eenum
\eprop 

\bpf 
\hfill
\benum
\item
The claim is a special case of \S\ref{sch-case} \ref{l-c-H}.
 
 \item
The claim is a special case of \Cref{sm-open} \ref{sm-open-b}. 
 
 \item
 If $G_{R/\fm}$ is a scheme, then it is separated, so $\Delta_{\bbB G_{R/\fm}}$ is also separated by \Cref{Delta-BG}~\ref{Delta-BG-b}. The local compactness of $H^1(R, G)$ therefore follows from \Cref{sm-open} \ref{sm-open-b}. 
 \qedhere
 \eenum
\epf

We proceed to record conditions that ensure Hausdorffness of $H^1(R, G)$. For conditions that ensure Hausdorffness of $H^0(R, G)$, see \Cref{H-supreme}.

\bprop\lab{H1-lch}
Suppose that $R$ is Henselian and \'{e}tale-open. 

\begin{enumerate}[label={(\arabic*)}] 
\item \lab{H1-lch-1}
If $G$ is smooth, or

\item \lab{H1-lch-2}
If $R$ is finite-closed (in the sense of \S\ref{p-c-f-c}) and $G$ is finite, or

\item \lab{H1-lch-3}
If $R$ is finite-closed, $R^\times$ is closed in $R$, and $G$ is quasi-finite and separated, or

\item \lab{H1-lch-4}
If $R$ is a local field and $G$ is proper,~or

\item \lab{H1-lch-5}
If $R$ is the ring of integers of a nonarchimedean local field and $G$ is quasi-compact and separated, or

\item \lab{H1-lch-6}
If $R$ is proper-closed, $R^\times$ is closed in $R$, and $G$ is quasi-affine,
\eenum
then $H^1(R, G)$ is Hausdorff.
\eprop

\bpf 
In case \ref{H1-lch-1}, the claim follows from \Cref{sm-discrete} \ref{disc-1}. In cases \ref{H1-lch-2}--\ref{H1-lch-6}, the claim follows from \Cref{H-supreme}~\ref{H-5}--\ref{H-10} and \ref{H-8} (one also uses \Cref{Delta-BG}~\ref{Delta-BG-b}).
\epf

The following further Hausdorffness result was pointed out to us by Laurent Moret-Bailly. It proves in particular that $H^1(k, G)$ is Hausdorff if $k$ is a nonarchimedean local field and $G$ is a commutative $k$-group scheme of finite type.

\bprop\lab{MB-tip}
Let $k$ be the fraction field of a Henselian valuation ring, $\wh{k}$ the completion of $k$, and $G$ a $k$-group scheme of finite type. Suppose that the field extension $\wh{k}/k$ is separable. If the identity component $G^\circ$ is commutative or if $G_{\ov{k}}$ has no nonzero subtorus, then $H^1(k, G)$ is $\mathrm{T}_1$, i.e., its points are closed. In particular, if $G$ is commutative, then $H^1(k, G)$ is Hausdorff.
\eprop

\bpf
If $G$ is commutative, then $H^1(k, G)$ is an abelian topological group by \Cref{top-gp}~\ref{top-gp-H1} (together with \S\ref{et-op}~\ref{ff-Hens-val}), so the last sentence follows from the rest. 

If $k$ is discrete, then so is $H^1(k, G)$. If $k$ is nondiscrete, then \cite{GGMB14}*{Thm.~1.2~(1)~(c)} proves that for every finite type $k$-scheme $Y$ and every $G_Y$-torsor $f\colon X \ra Y$ the image of the induced map 
\[
f(k)\colon X(k) \ra Y(k)
\]
is closed. In other words, for every $k$-morphism $Y \xra{f} \bbB G$, the preimage of the neutral class of $\bbB G(k)$ is closed in $Y(k)$, i.e., the neutral element of $H^1(k, G)$ is closed. The closedness of the other elements of $H^1(k, G)$ then follows from \Cref{twist} because the assumptions on $G^\circ$ or on $G_{\ov{k}}$ may be checked after replacing $G$ by its base change to a finite subextension of $\ov{k}/k$ and are therefore inherited by the inner forms ${}_TG$.
\epf

We close \S\ref{topo-coho} with a result that describes topological properties of $H^1(R, G)$ in the case when $R$ is the ring of integers of a nonarchimedean local field.

\bprop \lab{pullb-o-k}
Let $k$ be a nonarchimedean local field, $\fo$ its ring of integers, $\bF$ its residue field, and $G$ a flat and locally of finite type $\fo$-group algebraic space. The pullback map 
\[
p\colon H^1(\fo, G) \ra H^1(k, G)
\]
is continuous. If $G_\bF$ is a scheme, then $p$ is open. If, in addition to $G_\bF$ being a scheme, $G$ is quasi-compact (resp.,~quasi-compact and separated), then $H^1(\fo, G)$ is compact (resp.,~compact and Hausdorff).\eprop

\bpf
The continuity of $p$ follows from the aspect \ref{R-Rpr-cts} of \Cref{basic} applied to $X = \bbB G$. If $G_\bF$ is a scheme, then $\Delta_{\bbB G_{\bF}}$ is separated by \Cref{Delta-BG}~\ref{Delta-BG-b}, and hence $p$ is open by \Cref{sm-open}~\ref{sm-open-c}. If $G_\bF$ is a scheme and $G$ is quasi-compact, then $G_\bF$ is a group scheme of finite type over a finite field, so $H^1(\bF, G)$ is finite and the compactness of $H^1(\fo, G)$ follows from \Cref{sm-open} \ref{compact}. The Hausdorffness claim follows from \Cref{H1-lch} \ref{H1-lch-5}.
\epf


\section{Topological properties of maps in cohomology sequences} \lab{cont-les}

\bpp[The setup] \lab{setup}
As in \S\ref{topo-coho}, let $R$ be a local topological ring that satisfies \ref{alpha}--\ref{beta} of \S\ref{R-def}. Similarly to \Cref{sm-supreme}, let 
\[
\iota\colon H \hra G
\]
be a monomorphism of flat $R$-group algebraic spaces locally of finite presentation, use $\iota$ to identify $H$ with a subspace of $G$, and let 
\[
Q \ce G/H
\]
be the resulting quotient $R$-algebraic space, which is flat and locally of finite presentation by \Cref{sm-supreme} \ref{sm-supreme-a} and \cite{SP}*{\href{http://stacks.math.columbia.edu/tag/06ET}{06ET} and \href{http://stacks.math.columbia.edu/tag/06EV}{06EV}}. 

In this section, we investigate the topological properties of the maps in the long exact sequence
\be \lab{les} \tag{$\bigstar$}
1 \ra H(R) \xra{\iota(R)} G(R) \ra Q(R) \ra H^1(R, H) \ra H^1(R, G) \xra{x} H^1(R, Q) \xra{y} H^2(R, H) \xra{z} \dotsc,
\ee
where $x$ (resp.,~$y$) is defined if $\iota(H)$ is normal (resp., central) in $G$, and $z$ and the subsequent maps are defined if $G$ is commutative. 

As in \S\ref{topo-coho}, we have nothing to add to \Cref{sm-discrete} \ref{disc-2}--\ref{disc-4} concerning the topological study of $H^n(R, -)$ with $n \ge 2$. Throughout \S\ref{cont-les} we therefore mostly omit such $n$ from consideration.
\epp

\bprop \lab{les-cont}
All the maps in \eqref{les} are continuous (whenever defined). 
\eprop

\bpf
All the maps have underlying morphisms of presheaves, so the discussion of \S\ref{presheaf-case} applies. 
\epf

In \Cref{open-supreme} (resp.,~\Cref{cl-supreme}), we record conditions which ensure that the maps appearing in the long exact sequence \eqref{les} are open (resp.,~closed). Knowledge of these conditions will be handy in the proofs of the comparison with the \v{C}ech topology carried out in \S\S\ref{comparison}--\ref{comparison-2}.

\bprop \lab{open-supreme}
Suppose that $(R, \fm)$ is Henselian and \'{e}tale-open. 
\benum
\item \lab{open-supreme-a}
If $H \ra \Spec R$ is smooth, then $G(R) \ra Q(R)$ is open.

\item \lab{open-supreme-b}
If $G \ra \Spec R$ is smooth, then $Q(R) \ra H^1(R, H)$ is open.

\item \lab{open-supreme-c}
If $Q \ra \Spec R$ is smooth and $H_{R/\fm}$ is a scheme, then $H^1(R, H) \ra H^1(R, G)$ is open.

\item \lab{open-supreme-d}
If $\iota(H)$ is normal in $G$ and $G_{R/\fm}$ is a scheme, then $H^1(R, G) \xra{x} H^1(R, Q)$ is open.
\eenum
\eprop

\bpf
The claims follow by combining \Cref{sm-supreme} for $\cP = \text{`smooth'}$ and \Cref{sm-open}~\ref{sm-open-a} (together with the last sentence of \Cref{main-case} to ensure that the diagonal assumption of \Cref{sm-open} is met).
\epf

\bprop \lab{cl-supreme}
Consider the conditions: 
\benuma
\item \lab{Pn-c-cond}
$R$ is a local field;
\item \lab{Pn-2-cond}
$R$ is the ring of integers of a nonarchimedean local field;
\item \lab{f-c-cond}
$R$ is finite-closed;
\item \lab{p-c-cond}
$R$ is proper-closed;
\item \lab{f-c-cond-2}
$R$ is finite-closed and $R^\times$ is closed in $R$;
\item \lab{p-c-cond-2}
$R$ is proper-closed and $R^\times$ is closed in $R$. 
\eenum

Then
\benum
\item \lab{cl-supreme-a}
If $H \ra \Spec R$ is proper (resp.,~quasi-compact and separated, resp.,~finite, resp.,~proper and $G$ is a scheme, resp.,~quasi-finite and $G$ is a scheme, resp.,~quasi-compact and $G$ is a scheme) and $R$ satisfies \ref{Pn-c-cond} (resp.,~\ref{Pn-2-cond}, resp.,~\ref{f-c-cond}, resp.,~\ref{p-c-cond}, resp.,~\ref{f-c-cond-2}, resp.,~\ref{p-c-cond-2}), then 
\[
G(R) \ra Q(R) \qq \text{is closed;}
\]

\item \lab{cl-supreme-b}
If $G \ra \Spec R$ is proper (resp.,~quasi-compact and separated, resp.,~finite, resp.,~proper, $Q$ is a scheme, and $H$ is a quasi-affine scheme, resp.,~quasi-finite, $Q$ is a scheme, and $H$ is a quasi-affine scheme, resp.,~quasi-compact, $Q$ is a scheme, and $H$ is a quasi-affine scheme) and $R$ satisfies \ref{Pn-c-cond} (resp.,~\ref{Pn-2-cond}, resp.,~\ref{f-c-cond}, resp.,~\ref{p-c-cond}, resp.,~\ref{f-c-cond-2}, resp.,~\ref{p-c-cond-2}), then 
\[
Q(R) \ra H^1(R, H) \qq \text{is closed;}
\]

\item \lab{cl-supreme-c}
If $Q \ra \Spec R$ is proper (resp.,~quasi-compact and separated, resp.,~finite, resp.,~quasi-finite and separated, resp.,~quasi-affine) and $R$ satisfies \ref{Pn-c-cond} (resp.,~\ref{Pn-2-cond}, resp.,~\ref{f-c-cond}, resp.,~\ref{f-c-cond-2}, resp.,~\ref{p-c-cond-2}), then 
\[
H^1(R, H) \ra H^1(R, G) \qq \text{is closed.}
\]
\eenum
\eprop

\bpf
Due to \Cref{sm-supreme}, depending on whether we are in the situation \ref{cl-supreme-a}, or \ref{cl-supreme-b}, or \ref{cl-supreme-c}, the morphism $G \ra Q$, or $Q \ra \bbB H$, or $\bbB H \ra \bbB G$ inherits the geometric properties imposed on $H$, or $G$, or $Q$, respectively. Therefore,

\q\ In the case \ref{Pn-c-cond}, the claims \ref{cl-supreme-a}, \ref{cl-supreme-b}, and \ref{cl-supreme-c} follow from \Cref{Pn-cpct} \ref{Pn-cpct-d};

\q\ In the case \ref{Pn-2-cond}, the claims \ref{cl-supreme-a}, \ref{cl-supreme-b}, and \ref{cl-supreme-c} follow from \Cref{turbo-c};

\q\ In the case \ref{f-c-cond}, the claims \ref{cl-supreme-a}, \ref{cl-supreme-b}, and \ref{cl-supreme-c} follow from \Cref{prop-cl} \ref{prop-cl-b};

\q\ In the case \ref{p-c-cond}, the claims \ref{cl-supreme-a} and \ref{cl-supreme-b} follow from \Cref{prop-cl} \ref{prop-cl-a} (for the claim \ref{cl-supreme-b}, one uses the quasi-affineness of $H$ through \Cref{Delta-BG} \ref{Delta-BG-b} to argue that $\Delta_{\bbB H/R}$ is quasi-affine, so that $Q \ra \bbB H$ is representable by schemes);

\q\ In the case \ref{f-c-cond-2}, the claims \ref{cl-supreme-a}, \ref{cl-supreme-b}, and \ref{cl-supreme-c} follow from \Cref{gamma} \ref{gamma-b}, which applies because finite-closedness implies Hausdorffness (see \S\ref{p-c-f-c}) whereas the representability by schemes of the morphism in question is argued for the claim \ref{cl-supreme-b} as in the proof of the case \ref{p-c-cond}, and is argued for the claim \ref{cl-supreme-c} via the fact \cite{LMB00}*{A.2} that a quasi-finite separated algebraic space over a scheme is a scheme;

\q\ In the case \ref{p-c-cond-2}, the claims \ref{cl-supreme-a}, \ref{cl-supreme-b}, and \ref{cl-supreme-c} follow from \Cref{gamma} \ref{gamma-c}, which applies thanks to the reasoning analogous to that of the proof of the case \ref{f-c-cond-2}.
\epf

\brems
\remi \lab{quot-top}
Consider a single statement of \Cref{open-supreme} or \Cref{cl-supreme} and let $t$ denote the map of \eqref{les} studied there. A continuous open (resp.,~closed) surjection is a quotient map, so if in addition to the assumptions of the statement in question also the map of \eqref{les} succeeding $t$ is defined and vanishes, then the target of $t$ inherits a quotient topology from the source. 

\remi \lab{BT-approach}
The approach of \cite{BT14} is based on an analogue of Remark \ref{quot-top} as follows. Op.~cit.~assumes that $R$ is the fraction field of a complete rank $1$ valuation ring and $H$ is affine, embeds $H$ into a suitable smooth $G$ for which $H^1(R, G)$ vanishes (e.g.,~into $GL_n$), and endows $H^1(R, H)$ with the resulting quotient topology. It then argues in \cite{BT14}*{2.1.1} that the choice of $G$ did not matter and claims in \cite{BT14}*{3.1} that the quotient topology agrees with the \v{C}ech topology (defined in \cite{BT14}*{2.2}). Conclusions about the \v{C}ech topology use this analysis---an analogue of \Cref{les-cont} is \cite{BT14}*{5.1}, that of \Cref{open-supreme} \ref{open-supreme-c} is \cite{BT14}*{5.1.2}, and that of \Cref{sm-discrete} \ref{disc-1} is \cite{BT14}*{5.1.3}. However, as mentioned in Remark \ref{BT-rem}, the \v{C}ech topology used in \cite{BT14} is a priori different from the one used here and in \cite{Mil06}*{III.\S6}, so we cannot use the results cited above to shorten the proofs given in \S\S\ref{comparison}--\ref{comparison-2} below.

\remi
In contrast to \Cref{open-supreme} \ref{open-supreme-d}, $x$ need not be closed even if $R$ is a local field and $H$, $G$, and $Q$ are commutative and finite. For example, consider 
\[
\qqq 0 \ra \gA_p \ra \gA_{p^2} \ra \gA_p \ra 0
\]
over $R \ce \bF_p((t))$. The vanishing of $H^n(R, \bG_a)$ for $n \ge 1$ combines with Remark \ref{quot-top} to show that 
\[
\qqq H^1(R, \gA_{p^2}) = R/R^{p^2} \qq \text{and} \qq H^1(R, \gA_p) = R/R^p
\]
both algebraically and topologically, and, moreover, that $x$ identifies with the quotient map 
\[
\qqq R/R^{p^2} \surjects R/R^p.
\]
The latter is not closed: 
\[
\qqq \{ t^n + t^{-pn}\}_{n \ge 1\text{ and } p\nmid n} \subset R
\]
has a closed image in $R/R^{p^2}$ but not in $R/R^p$. Also, the $p$-power map $R/R^p \ra R/R^{p^2}$ is not open, so the smoothness assumption of \Cref{open-supreme} \ref{open-supreme-c} cannot be dropped. Neither can that of \Cref{open-supreme} \ref{open-supreme-b} because $R/R^p$ is nondiscrete. 
\erems


\section{Comparison with the \v{C}ech topology on $H^1(k, G)$} \lab{comparison}

The main goal of the present section is \Cref{main-Cech}, which shows that the topology on $H^1(k, G)$ introduced in \S\ref{H1-top} agrees with the \v{C}ech topology used by other authors. 

\bpp[The \v{C}ech topology] \lab{Cech-topo}
Let $k$ be a Hausdorff topological field and $G$ a locally of finite type $k$-group scheme. For a finite subextension $\ov{k}/L/k$ and $n \ge 0$, set 
\[
L_n \ce \tensor_{i = 0}^n L \qq \text{(tensor product over $k$).}
\]
Consider the restriction of scalars 
\[
C^n \ce \Res_{L_n / k} (G_{L_n}).
\]
By \cite{EGA IV3}*{9.1.5} and \cite{CGP10}*{A.3.5}, the connected components of $G$ are quasi-projective. Therefore, every finite set of points of $G_{L_n}$ is contained in an open affine, and hence $C^n$ is a locally of finite type $k$-group scheme by \cite{BLR90}*{7.6/4}. We call $C^n$ the \emph{scheme of $n$-cochains} (for $G$ with respect to $L/k$). The coboundaries 
\[
d^n\colon C^n \ra C^{n + 1},
\]
which are defined using the usual formulas, are morphisms of $k$-schemes; they are even group homomorphisms if $G$ is commutative. Let $e_n\colon \Spec k \ra C^{n}$ be the unit section and 
\[
Z^n \ce C^n \times_{d^n, C^{n+1}, e_{n + 1}} \Spec k 
\]
the \emph{scheme of $n$-cocycles} (for $G$ with respect to $L/k$).
 
The group $C^0(k)$ acts on the right on $Z^1(k)$ by 
\[
(z^1, c^0) \mapsto p_1^*(k)(c^0)\i z^1 p_0^*(k)(c^0)
\]
where $p_i^*\colon C^0 \ra~C^1$ is the map induced by the $i\th$ projection $\Spec L_1 \ra \Spec L$ (so $d^0 = (p_1^*)\i p_0^*$). We endow the orbit space 
\[
H^1(L/k, G) \ce Z^1(k)/C^0(k)
\]
with the quotient topology. Its points correspond to $G$-torsors trivialized by $L/k$, so 
\[
H^1(L/k, G) \hra H^1(k, G).
\]
If $L\pr$ is contained in $L$, then the inclusion $H^1(L\pr/k, G) \hra H^1(L/k, G)$ is continuous because it is induced by a $k$-morphism $Z^1_{L\pr/k} \ra Z^1_{L/k}$. The transition maps in 
\[
\textstyle{H^1(k, G) = \varinjlim_{L} H^1(L/k, G)}
\]
are therefore continuous, and we write $H^1(k, G)_{\text{\v{C}ech}}$ for $H^1(k, G)$ endowed with the direct limit topology. A subset $U \subset H^1(k, G)_{\text{\v{C}ech}}$ is open if and only if so is its preimage in every $Z^1(k)$.

If $G$ is commutative, then $(C^n(k), d^n(k))_{n \ge 0}$ is a complex of abelian groups. We endow its cohomology groups $H^n(L/k, G)$ with the subquotient topology. If $G$ is locally of finite type, then 
\[
\textstyle{H^n(k, G) = \varinjlim_L H^n(L/k, G)}
\]
by \cite{Sha72}*{p.~208, Thm.~42},\footnote{
Loc.~cit.~assumes that $G$ is of finite type but the method of proof there continues to work for locally of finite type $G$ because such a $G$ is an extension of a smooth $k$-group by a finite $k$-group due to \Cref{SGA3-input}.
} 
and we write $H^n(k, G)_{\text{\v{C}ech}}$ for $H^n(k, G)$ endowed with the direct limit topology. 

The agreement of our definition of the \v{C}ech topology with the one used in \cite{Mil06}*{III.\S6} follows from \Cref{weil-weil} \ref{weil-weil-a}. The absence of the \v{C}ech subscript indicates the topologies of \S\S\ref{H1-top}--\ref{Hi-top}.
\epp

\brem \lab{H0-homeo}
The bijection 
\[
H^0(k, G)_{\text{\v{C}ech}} \ra H^0(k, G)
\]
is a homeomorphism. In fact, both topologies agree with the topology on $G(k)$: for $H^0(k, G)_{\text{\v{C}ech}}$, this follows from \S\ref{sch-case} \ref{R-Rpr-c}; for $H^0(k, G)$, this follows from the definition given in \S\ref{H1-top}.
\erem

\blem \lab{Cech-cont}
The bijection 
\[
b_G\colon H^1(k, G)_{\text{\upshape{\v{C}ech}}} \ra H^1(k, G)
\]
is continuous.
\elem

\bpf
We need to argue that for every $L$ the composite 
\[
Z^1(k) \ra H^1(k, G)_{\text{\v{C}ech}} \ra H^1(k, G)
\]
is continuous (as before, $Z^1$ depends on $L$). For this, it suffices to exhibit an underlying morphism 
\[
Z^1 \xra{a} \bbB G.
\]
On the category of $k$-schemes, $Z^1$ represents the functor of $1$-cocycles for $G$ with respect to $L/k$. The universal $1$-cocycle corresponding to $\id_{Z^1}$ gives rise to a torsor under $G_{Z^1} \ra Z^1$, and this torsor corresponds to $a$.
\epf

\begin{cor} \lab{sm-homeo}
If $k$ is \'{e}tale-open and $G \ra \Spec k$ is smooth, then $b_G$ is a homeomorphism and $H^1(k, G)_{\text{\upshape{\v{C}ech}}}$ is discrete.
\end{cor}

\bpf
By \Cref{sm-discrete} \ref{disc-1}, $H^1(k, G)$ is discrete, so continuity of $b_G$ gives the claim.
\epf


\blem \lab{H1-cont}
For a homomorphism $f\colon H \ra G$ of locally of finite type $k$-group schemes, the resulting morphism 
\[
H^1(k, H)_{\text{\upshape{\v{C}ech}}} \ra H^1(k, G)_{\text{\upshape{\v{C}ech}}}
\]
is continuous.
\elem

\bpf
For a finite subextension $\ov{k}/L/k$, write $f_{Z^1}\colon Z^1_H \ra Z^1_G$ for the $k$-morphism of schemes of $1$-cocycles induced by $f$. The claim results from the commutativity of
\[
\xymatrix{
Z^1_H(k) \ar[r]^{f_{Z^1}(k)} \ar[d] & Z^1_G(k) \ar[d] \\
H^1(k, H)_{\text{\v{C}ech}} \ar[r] & H^1(k, G)_{\text{\v{C}ech}}
}
\]
and the continuity of $f_{Z^1}(k)$.
\epf

\blem \lab{twist-Cech}
Let $T$ be a right $G$-torsor and let ${}_TG \ce \Aut_G(T)$ be the corresponding inner form of $G$. Twisting by $T$ induces a homeomorphism 
\[
H^1(k, G)_{\text{\upshape{\v{C}ech}}} \cong H^1(k, {}_TG)_{\text{\upshape{\v{C}ech}}}
\]
that sends the class of $T$ to the neutral class. 
\elem

\bpf
We write $Z^1_{L/k, G}$ for the $Z^1$ of \S\ref{Cech-topo} and choose a finite subextension $\ov{k}/L/k$ that trivializes $T$ and an $x \in Z^1_{L/k, G}(k)$ that gives the class of $T$ in $H^1(L/k, G)$. For every finite subextension $\ov{k}/L\pr/L$, right multiplication by the pullback of $x\i$ induces an isomorphism 
\[
i\colon Z^1_{L\pr/k, G} \isomto Z^1_{L\pr/k, {}_TG}
\]
(compare with \cite{Ser02}*{I.\S5, Prop.~35~bis}). Moreover, the homeomorphism $i(k)$ intertwines the actions of $C^0_{L\pr/k, G}(k)$ and $C^0_{L\pr/k, {}_TG}(k)$. Thus, the resulting ``twist by $T$'' bijection 
\[
H^1(L\pr/k, G) \isomto H^1(L\pr/k, {}_TG)
\]
is a homeomorphism. It remains to take $\varinjlim_{L\pr}$. 
\epf

\blem \lab{fin-con}
Let $K$ be a field of characteristic $p$ and let $G$ be a finite connected $K$-group scheme. 
\benum
\item \lab{fin-con-a}
Every $G$-torsor $X \ra \Spec K$ for the \'{e}tale topology is trivial.

\item \lab{fin-con-b}
For some $n \ge 0$ with $p^n \le (\#G)!$, every $G$-torsor $X \ra \Spec K$ for the fppf topology is trivial over $K^{1/p^n} \subset \ov{K}$ (with the convention that $0^0 = 1$). In particular, if $K$ is perfect, then $H^1(K, G)$ is trivial.
\eenum
\elem

\bpf
If $p = 0$, then $G$ is trivial, so we assume that $p > 0$.
\benum
\item
If $A$ is a reduced $K$-algebra, then $G(A)$ is a singleton. We take $A = L \tensor_K L$ for a finite separable extension $L/K$ to see that there are no nontrivial $G$-torsors trivialized by $L/K$.

\item
For a subextension $\ov{K}/K\pr/K$ of degree at most $\#G$, let $N \subset \ov{K}$ be its normal closure, so $[N : K] \le (\#G)!$.
Let $n$ be such that $p^n \le (\#G)!$ but $p^{n + 1} > (\#G)!$. Set $L \ce K^{1/p^n}$.

If $K\pr$ is a residue field of $X$, then, by \cite{BouA}*{V.73, Prop.~13}, $NL/L$ is separable. Thus, $X_L$ trivializes over an \'{e}tale cover of $\Spec L$. Since $G_L$ is connected, it remains to apply \ref{fin-con-a}. \qedhere
\eenum
\epf

\blem \lab{SGA3-input}
Let $K$ be a field. Every locally of finite type $K$-group scheme $G$ is an extension
\[
1 \ra H \ra G \ra Q \ra 1
\]
with $Q$ smooth and $H$ finite, connected, and closed and normal in $G$.
\elem

\bpf
If $\Char K = 0$, then $G$ is smooth by Cartier's theorem \cite{SGA3Inew}*{VI$_{\text{\upshape{B}}}$, 1.6.1}, so the trivial subgroup $H = 1$ suffices. If $\Char K > 0$, then \cite{SGA3Inew}*{VII$_{\text{\upshape{A}}}$, 8.3} gives the claim (with $H$ a large enough relative Frobenius kernel).
\epf

\bpp[The ground field] \lab{ground}
For the rest of \S\ref{comparison}, $k$ is a Hausdorff topological field that is \'{e}tale-open and finite-closed in the sense of \S\ref{et-op} and \S\ref{p-c-f-c} and such that every finite extension $L$ of $k$ is finite-closed when endowed with the topology resulting from some (equivalently, any) $k$-module identification $L \cong k^n$. For instance, as justified in \S\ref{et-op} and \S\ref{p-c-f-c}, $k$ could be a local field or a field that is complete with respect to a nonarchimedean valuation of rank $1$. We set $p \ce \Char k$.
\epp

\blem \lab{BT-fix}
Assume that $[k : k^p] < \infty$ if $p > 0$. For a locally of finite type $k$-group scheme $G$, a finite connected closed subgroup $H \le G$ (that need not be normal), and $Q \ce G/H$ (which is a locally of finite type $k$-scheme by \cite{SGA3Inew}*{VI$_{\text{\upshape{A}}}$, 3.2}), the connecting map 
\[
\delta(k)\colon Q(k) \ra H^1(k, H)_{\text{\upshape{\v{C}ech}}}
\]
is continuous.
\elem

\bpf
We use \Cref{fin-con} \ref{fin-con-b} to fix a finite subextension $\ov{k}/L/k$ for which $H^1(L/k, H) \hra H^1(k, H)$ is bijective. We indicate the relevant groups by subscripts and consider the Cartesian square
\[
\xymatrix{
F \ar@{^(->}[r]^-{i}\ar[d]_-{g} & C^0_G \ar[d]^-{d^0_G} \\
Z^1_H \ar@{^(->}[r] & Z^1_G,
}
\]
in which the horizontal maps are closed immersions thanks to \cite{BLR90}*{7.6/2 (ii)}.

\bcl \lab{BT-fix-1}
There is a Cartesian square (for the bottom map, recall that $C^0_G = \Res_{L/k} (G_L)$)
\be\ba \lab{F-cart-2}
\xymatrix{
F \ar@{^(->}[d]_-{i}\ar[r]^-{f} & Q \ar@{^(->}[d] \\
C^0_G \ar[r] & \Res_{L/k}(Q_L).
}
\ea\ee
\ecl

\bpf
The claim is a manifestation of fppf descent: for a $k$-scheme $S$, the two pullbacks 
\[
p_0^*(x), p_1^*(x) \in C^1_G(S) = G(S_L \times_S S_L)\quad \text{ of an element }\quad x\in C^0_G(S) = G(S_L)
\]
satisfy
\[
p_0^*(x) = p_1^*(x) y \quad \text{ for some }\quad y \in Z^1_H(S) \subset H(S_L \times_S S_L)
\]
if and only if the image of $x$ in $Q(S_L)$ lands in $Q(S)$.
\epf

\bcl \lab{BT-fix-2}
$f(k)\colon F(k) \ra Q(k)$ is continuous, closed, and surjective.
\ecl

\bpf
Continuity is clear (from \S\ref{sch-case} \ref{cont}). By the assumptions of \S\ref{ground}, $L$ is finite-closed, so 
\[
G(L) \ra Q(L)
\]
is closed. To deduce the same for $f(k)$, we evaluate \eqref{F-cart-2} on $k$-points and apply \S\ref{sch-case}~\ref{c-i-H} to $i$, \Cref{weil-weil}~\ref{weil-weil-a} to $C^0_G$, and \S\ref{sch-case} \ref{R-Rpr-cts} to $Q(k) \ra Q(L)$. By the choice of $L$, every $k$-point of $Q$ lifts to an $L$-point of $G$, so \eqref{F-cart-2} also gives the surjectivity of $f(k)$.
\epf

\bcl \lab{BT-fix-3}
The following diagram commutes:
\be\ba\lab{BT-comm}
\xymatrix{
F(k) \ar[r]^-{f(k)} \ar[d]_-{g(k)} & Q(k) \ar[d]^-{\delta(k)} \\
Z^1_H(k) \ar[r]^-{\pi} & H^1(k, H)_{\text{\v{C}ech}}.
}
\ea\ee
\ecl

\bpf
We fix a $q \in Q(k)$, let $X \ra \Spec k$ be the corresponding fiber of $G \ra Q$, and view $X$ as an $H$-torsor trivialized by $L/k$. To build a 
\[
z \in Z^1_H(k) \subset H(L\tensor_k L)
\]
that gives the class of $X$ in $H^1(L/k, H)$, one takes any 
\[
x \in X(L) \overset{\eqref{F-cart-2}}{\subset} F(k) \subset G(L)
\]
and sets $z \ce d^0_G(x)$ to be the $H(L\tensor_k L)$-difference of the two pullbacks of $x$ to $X(L \tensor_k L)$. This gives the claim because $d_G^0(x) = g(k)(x)$. 
\epf

To prove the continuity of $\delta(k)$, let $Z \subset H^1(k, H)_{\text{\v{C}ech}}$ be closed. The continuity of $\pi$ and $g(k)$ gives the closedness of 
\[
g(k)\i(\pi\i(Z)) \overset{\eqref{BT-comm}}{=} f(k)\i(\delta(k)\i(Z)),
\]
so \Cref{BT-fix-2} gives that of $\delta(k)\i(Z)$.
\epf

\bthm \lab{main-Cech}
Let $k$ be a Hausdorff topological field of characteristic $p \ge 0$ such that $k$ is \'{e}tale-open, every finite extension of $k$ endowed with the $k$-vector space topology is finite-closed, and $[k : k^p] < \infty$ if $p > 0$. (For instance, $k$ could be a local field.) For a locally of finite type $k$-group scheme $G$, the bijection 
\[
b_G\colon H^1(k, G)_{\text{\upshape{\v{C}ech}}} \ra H^1(k, G)
\]
is a homeomorphism.
\ethm

\bpf
We begin with the case of a finite connected $G$. We embed $G$ into some $\GL_n$, and set $Q \ce \GL_n/G$. By \Cref{open-supreme} \ref{open-supreme-b} (with Remark \ref{quot-top}), the connecting map 
\[
\delta(k)\colon Q(k) \surjects H^1(k, G)
\]
is a quotient map. On the other hand, by \Cref{BT-fix}, 
\[
\delta(k)\colon Q(k) \surjects H^1(k, G)_{\text{\v{C}ech}}
\]
is continuous. Therefore, 
\[
b_G\colon H^1(k, G)_{\text{\v{C}ech}} \ra H^1(k, G)
\]
is open, and hence is also a homeomorphism thanks to its continuity supplied by \Cref{Cech-cont}.

We now turn to the general case. By \Cref{SGA3-input}, there is an exact sequence
\[
1 \ra H \ra G \ra Q \ra 1\quad\quad \text{with $H$ finite connected, $Q$ smooth, and $H$ closed and normal in $G$.}
\]
Its cohomology sequence gives the commutative diagram
\[
\xymatrix{
H^1(k, H)_{\text{\upshape{\v{C}ech}}} \ar[r]^-{\check{a}}\ar[d]^{b_H} & H^1(k, G)_{\text{\upshape{\v{C}ech}}} \ar[r]^-{\check{c}}\ar[d]^{b_G} &H^1(k, Q)_{\text{\upshape{\v{C}ech}}} \ar[d]^{b_Q} \\
 H^1(k, H) \ar[r]^{a} & H^1(k, G) \ar[r] &H^1(k, Q),
}
\]
in which the horizontal maps are continuous by \Cref{les-cont} and \Cref{H1-cont}.
By \Cref{sm-homeo}, $H^1(k, Q)_{\text{\upshape{\v{C}ech}}}$ is discrete, so the fibers of $\check{c}$ are open. By \Cref{Cech-cont}, $b_G$ is continuous, so it remains to argue that the restrictions of $b_G$ to the fibers of $\check{c}$ are open morphisms. In the case of the neutral fiber, the restriction $b_G|_{\im{\check{a}}}$ is open because $b_H$ is a homeomorphism by the finite connected case and $a$ is open by \Cref{open-supreme}~\ref{open-supreme-c}. We will argue that the openness of the restrictions of $b_G$ to the other fibers of $\check{c}$ follows from the settled case of the neutral fiber.

Let $T$ be a $G$-torsor, and let ${}_TG \ce \Aut_G(T)$ be the resulting inner form of $G$. By \Cref{twist} and \Cref{twist-Cech}, twisting by $T$ induces homeomorphisms 
\[
H^1(k, G)_{\text{\upshape{\v{C}ech}}} \cong H^1(k, {}_TG)_{\text{\upshape{\v{C}ech}}}
\]
and 
\[
H^1(k, G) \cong H^1(k, {}_TG),
\]
which, moreover, are compatible with $b_G$ and $b_{{}_TG}$. By \cite{Gir71}*{III.3.3.5}, these homeomorphisms identify the $\check{c}$-fiber containing the class of $T$ with the neutral fiber of 
\[
H^1(k, {}_TG)_{\text{\upshape{\v{C}ech}}} \ra H^1(k, {}_{T/H}Q)_{\text{\upshape{\v{C}ech}}},
\]
where $T/H$ is regarded as a $Q$-torsor and ${}_{T/H}Q \ce \Aut_Q(T/H)$ is the resulting inner form of $Q$. This identification achieves the desired reduction to the case of the neutral fiber because ${}_{T/H}Q$ is smooth and $\Ker({}_TG \ra {}_{T/H}Q)$ is finite connected, as may be checked after base change to a finite extension trivializing $T$.
\epf

\beg \lab{BT-counter}
We show that the assumption on $[k : k^p]$ cannot be dropped in \Cref{main-Cech}. We take 
\[
k = \bF_p(x_1, x_2, \dotsc)((t)) \qq \text{and} \qq G = \gA_p.
\]
Since $k$ is a complete discrete valuation ring, it satisfies the assumptions of \Cref{main-Cech} with the exception of the requirement on $[k : k^p]$. By \Cref{open-supreme} \ref{open-supreme-b} and Remark \ref{quot-top}, $H^1(k, G)$ identifies with $k/k^p$ equipped with the quotient topology. In particular, 
\[
\{ x_nt^{np}\}_{n \ge 1} \subset H^1(k, G)
\]
has $0$ as a limit point, and hence is not closed. On the other hand, only finitely many $x_n$ have a $p\th$ root in a fixed finite subextension $\ov{k}/L/k$, so the preimage of 
\[
\{ x_nt^{np}\}_{n \ge 1} \subset H^1(k, G)_{\text{\v{C}ech}}
\]
in $H^1(L/k, G)$ is finite. Since $C^0(k) = \alpha_p(L)$ is finite, the further preimage in $Z^1(k)$ is also finite, and hence closed thanks to \Cref{H-supreme} \ref{H-1}. In conclusion, 
\[
\{ x_nt^{np}\}_{n \ge 1} \subset H^1(k, G)_{\text{\v{C}ech}}
\]
is closed, so $b_G$ is not a homeomorphism.
\eeg


\section{Comparison with the \v{C}ech topology on $H^n(k, G)$ for $n \ge 2$} \lab{comparison-2}

The goal of \S\ref{comparison-2} is \Cref{Hi-disc}: for commutative $G$ and $n \ge 2$, the \v{C}ech topology on $H^n(k, G)$ is simply the discrete topology, which by \Cref{sm-discrete} \ref{disc-2} agrees with the topology of~\S\ref{Hi-top}.

\bpp[The ground field] \lab{ground-2}
Throughout \S\ref{comparison-2}, $k$ is a Hausdorff topological field that is \'{e}tale-open in the sense of \S\ref{et-op}. For instance, $k$ could be a local field or the fraction field of a Henselian valuation ring.
\epp

\bpp[The \v{C}ech topology]
Focusing on commutative $G$, we study the \v{C}ech topology defined in~\S\ref{Cech-topo}.  
\epp

We begin with the analogue of \Cref{sm-homeo} in the commutative case:

\blem \lab{sm-disc}
If $G$ is a commutative smooth $k$-group scheme, then $H^n(k, G)_{\text{\upshape{\v{C}ech}}}$ is discrete for $n \ge 1$.
\elem

\bpf
Fix a finite subextension $\ov{k}/L/k$. By \Cref{sm-open} \ref{sm-open-a} and \Cref{Groth-2}, the fibers of 
\[
Z^n(k) \ra H^n(L/k, G)
\]
are open. Thus, $H^n(L/k, G)$ is discrete, and hence so is $H^n(k, G)_{\text{\upshape{\v{C}ech}}}$.
\epf

\blem \lab{sm-embed}
For a field $K$, every commutative locally of finite type $K$-group scheme $H$ is a closed subgroup of a commutative smooth $K$-group scheme $G$, which can be chosen to be quasi-compact if so is $H$.
\elem

\bpf
By \Cref{SGA3-input}, there is a short exact sequence
\[
0 \ra H_0 \ra H \ra H\pr \ra 0\quad\quad \text{with $H_0$ finite connected and $H\pr$ smooth.}
\]
Let $H_0^D$ be the Cartier dual of $H_0$, so that $H_0 = \mathcal{H}om(H_0^D, \bG_m)$. Embed $H_0$ into the smooth $k$-group $G_0 \ce \Res_{H_0^D/k} \bG_m$ via the functor $\mathcal{H}om(H_0^D, \bG_m) \hra \mathcal{M}or(H_0^D, \bG_m)$. The pushout 
\[
\xymatrix{
0 \ar[r] &H_0 \ar@{^(->}[d] \ar[r] &H \ar@{^(->}[d]^-{i} \ar[r] &H\pr \ar@{=}[d] \ar[r] &0 \\
0 \ar[r] & G_0 \ar[r] &G \ar[r] &H\pr \ar[r] &0
}
\]
defines a desired $G$, because $G$ inherits $k$-smoothness from $G_0$ and $H\pr$.
\epf

\bthm \lab{Hi-disc}
If $H$ is a commutative locally of finite type $k$-group scheme, then $H^n(k, H)_{\text{\upshape{\v{C}ech}}}$ is discrete for $n \ge 2$. 
\ethm

\bpf
When needed, we indicate the relevant groups and subextensions by subscripts.

Fix an $a \in H^n(k, H)_{\text{\upshape{\v{C}ech}}}$ and its preimage $b \in Z^n_{L/k, H}(k)$ for a sufficiently large $L$. Let $V \subset Z^n_{L/k, H}(k)$ be the preimage of an open $U \subset H^n(k, H)_{\text{\upshape{\v{C}ech}}}$. The preimage of $a + U$ is the open $b + V$. Thus, the translation by $a$ map is open, and we are reduced to proving the openness of $\{ 0\} \subset H^n(k, H)_{\text{\upshape{\v{C}ech}}}$.

Use \Cref{sm-embed} to embed $H$ into a commutative smooth $k$-group scheme $G$ and set $Q\ce G/H$, which is a commutative smooth $k$-group scheme by \cite{SGA3Inew}*{VI$_{\text{\upshape{A}}}$,~3.2}. Let 
\[
\delta\colon H^{n - 1}(k, Q)_{\text{\upshape{\v{C}ech}}} \ra H^n(k, H)_{\text{\upshape{\v{C}ech}}}
\]
be the connecting homomorphism. By \Cref{sm-disc}, for a fixed finite subextension $\ov{k}/L/k$,
\[
X \ce \Ker(Z^{n - 1}_Q(k) \ra H^{n - 1}(k, Q)_{\text{\upshape{\v{C}ech}}})\quad \text{ is a nonempty open in } \quad Z^{n - 1}_Q(k).
\] 
Let $Y \subset C^{n - 1}_G(k)$ be the preimage of $X\subset Z^{n - 1}_Q(k)$. By \cite{Mil06}*{proof of III.6.1}, the nonempty
\be\lab{Mil-input}
d^{n - 1}_G(k)(Y) \subset Z^{n}_H(k) \quad\quad \text{lies in the preimage of}\quad\quad \delta(0) \in H^n(k, H)_{\text{\upshape{\v{C}ech}}}.
\ee

\bcl \lab{LkW}
$d^{n - 1}_G(k)(Y)$ is open in $Z^n_H(k)$.
\ecl

\bpf
The proof is similar to that of \Cref{BT-fix}. Define the $k$-scheme $F$ to be the fiber product
\[
\xymatrix{
F \ar@{^(->}[d]\ar[r]^-{f} & Z^{n - 1}_Q \ar@{^(->}[d] \\
C^{n - 1}_{G} \ar[r] & C^{n - 1}_Q.
}
\]
Then $Y = f(k)\i(X)$, so the continuity of $f(k)$ guarantees the openness of $Y$ in $F(k)$. 

Left exactness of the restriction of scalars and the definition of $F$ ensure that the outer square of
\[
\xymatrix{
F \ar@{^(->}[d]\ar[r]^-{d_F} & Z^n_H \ar@{^(->}[r]\ar@{^(->}[d] &C^n_H \ar@{^(->}[d] \\
C^{n - 1}_G \ar[r]^-{d^{n - 1}_G} & Z^n_G \ar@{^(->}[r] &C^n_G
}
\]
is Cartesian. Since so is the right square, the left one is also. Thus, $d_F$ is smooth because so is $d_G^{n - 1}$ by \Cref{Groth-2}. Since 
\[
d^{n - 1}_G(k)(Y) = d_F(Y),
\]
it remains to combine the openness of $Y \subset F(k)$ with \Cref{sm-open}~\ref{sm-open-a}. 
\epf

Since $Z^{n}_H(k)$ is a topological group, \eqref{Mil-input} combines with \Cref{LkW} to prove that the preimage of $\{ 0\} \subset H^n(k, H)_{\text{\upshape{\v{C}ech}}}$ in $Z^n_H(k)$ is open. But $L$ was arbitrary, so $\{ 0\} \subset H^n(k, H)_{\text{\upshape{\v{C}ech}}}$ is open.
\epf


\appendix

\section{Generalities concerning stacks of torsors} \lab{BG-genl}

We gather generalities about the classifying stack $\bbB G$ that play an important role in the arguments of the main body of the paper. Specifically, we discuss algebraicity in \S\ref{BG-alg}, properties of the diagonal in \Cref{Delta-BG} (which we use without explicit notice), smoothness in \Cref{BG-sm}, and properties of morphisms arising from a short exact sequence in \Cref{sm-supreme}.

\bpp[Stacks of torsors] \lab{BG-alg}
Let $S$ be a scheme and let $G$ be an $S$-group algebraic space. Consider the classifying $S$-stack $\bbB G$, whose groupoid of $S\pr$-points for an $S$-scheme $S\pr$ is the category of right torsors under $G_{S\pr} \ra S\pr$ for the fppf topology (the torsors are algebraic spaces by \cite{SP}*{\href{http://stacks.math.columbia.edu/tag/04SK}{04SK}}, or already by \cite{LMB00}*{10.4.2} if $S$ and $G \ra S$ are quasi-separated). By \cite{SP}*{\href{http://stacks.math.columbia.edu/tag/06FI}{06FI}} (cf.~also \cite{LMB00}*{10.13.1}), $\bbB G$ is algebraic if $G \ra S$ is flat and locally of finite presentation. Let $S \xra{e} \bbB G$ be the morphism that corresponds to the trivial torsor. 

If $G$ is commutative, then the naturality of the contracted product and opposite torsor constructions shows that $\bbB G$ is a group object in the $2$-category of stacks over $S$. \epp

\blem \lab{Delta-BG}
Let $S$ be a scheme and let $G$ be an $S$-group algebraic space.
\benum
\item \lab{Delta-BG-a}
The following square is $2$-Cartesian:
\[
\xymatrix{
S \ar[rr]^{e}\ar[d]^{e} && \bbB G \ar[d]^{\id \times e} \\
\bbB G \ar[rr]^-{\Delta_{\bbB G/S}} && \bbB G \times_S \bbB G.
}
\]
\item \lab{Delta-BG-b}
Let $\cP$ be a property of morphisms of algebraic spaces that is stable under base change and is fppf local on the base. If $G \ra S$ has $\cP$, then so do $\Delta_{\bbB G/S}$ and $S \xra{e} \bbB G$.
\eenum

\elem

\bpf \hfill
\benum
\item
Inspect the definition of a $2$-fiber product.

\item
Checking that $\Delta_{\bbB G/S}$ has $\cP$ amounts to checking that for every $S$-scheme $S\pr$ and torsors $x$ and $y$ under $G_{S\pr} \ra S\pr$, the algebraic space $\Isom_{G_{S\pr}}(x, y) \ra S\pr$ parametrizing torsor isomorphisms has $\cP$. But $\Isom_{G_{S\pr}}(x, y) \ra S\pr$ is a right torsor for the fppf topology under an inner form of $G_{S\pr} \ra S\pr$ and hence inherits $\cP$. The claim that concerns $e$ now follows from \ref{Delta-BG-a}. \qedhere
\eenum
\epf

\bprop \lab{BG-sm}
Let $S$ be a scheme and let $G$ be a flat and locally of finite presentation $S$-group algebraic space. The algebraic stack $\bbB G$ is $S$-smooth.
\eprop

\bpf
By \Cref{Delta-BG} \ref{Delta-BG-b}, $S \xra{e} \bbB G$ is faithfully flat and locally of finite presentation, so the algebraic stack analogue of \cite{SP}*{\href{http://stacks.math.columbia.edu/tag/05B5}{05B5}} (or of \cite{EGAIV4}*{17.7.7}) applies to 
\[
S \xra{e} \bbB G \ra S. \qedhere
\]
\epf

In the following result, $\cP = \text{`smooth'}$ is the case of most interest to us.

\bprop \lab{sm-supreme}
Let $S$ be a scheme, $H$ a flat and locally of finite presentation $S$-group algebraic space, $\iota\colon H \hra G$ a monomorphism of $S$-group algebraic spaces, and  $Q \ce G/H$ the resulting quotient homogeneous space (which is an $S$-algebraic space by \cite{SP}*{\href{http://stacks.math.columbia.edu/tag/06PH}{06PH}} or, in many situations, already by \cite{LMB00}*{8.1.1 and 10.6}). Let $\cP$ be a property of morphisms of algebraic spaces that is stable under base change and is fppf local on the base.
\benum  
\item \lab{sm-supreme-a}
If $H \ra S$ has $\cP$, then $G \ra Q$ has $\cP$.

\item \lab{sm-supreme-b}
If $G \ra S$ has $\cP$, then the representable $Q \ra \bbB H$ that corresponds to $G \ra Q$ has $\cP$.

\item \lab{sm-supreme-c}
If $Q \ra S$ has $\cP$ and $G \ra S$ is flat and locally of finite presentation, then $\bbB H \xra{\bbB \iota} \bbB G$ is representable and has $\cP$. 

\item \lab{sm-supreme-d}
If $\iota(H)$ is normal in $G$ and $G \ra S$ is flat and locally of finite presentation, then $\bbB G \ra \bbB Q$ is smooth.
\eenum 
\eprop

\bpf 
For \ref{sm-supreme-a}, note that $G \ra Q$ is a torsor under $H \times_S Q \ra Q$ for the fppf topology. In \ref{sm-supreme-b}--\ref{sm-supreme-d}, the stacks in question are algebraic by the discussion in \S\ref{BG-alg} (for \ref{sm-supreme-d}, $Q$ is flat and locally of finite presentation by \ref{sm-supreme-a} and \cite{SP}*{\href{http://stacks.math.columbia.edu/tag/06ET}{06ET} and \href{http://stacks.math.columbia.edu/tag/06EV}{06EV}}). Moreover, by \Cref{Delta-BG} \ref{Delta-BG-b}, $e$ in  
\[
\xymatrix{
 G \ar[r]\ar[d] & Q \ar[d]  && Q \ar[r]\ar[d] & \bbB H \ar[d]^{\bbB \iota} & \ar@{}[d]^{\text{\normalsize{\ and}}} && \bbB H \ar[r]\ar[d] & \bbB G \ar[d] \\
 S \ar[r]^-{e} & \bbB H, && S \ar[r]^-{e} & \bbB G, &&& S \ar[r]^-{e} & \bbB Q
}
\]
is flat and locally of finite presentation in cases \ref{sm-supreme-b}, \ref{sm-supreme-c}, and \ref{sm-supreme-d}. Thus, to conclude using fppf descent (and \Cref{BG-sm} for \ref{sm-supreme-d}), we argue that the squares are $2$-Cartesian (in the respective~cases):
\benum \addtocounter{enumi}{1}
\item 
For a morphism $S\pr \xra{x} Q$ with $S\pr$ an $S$-scheme, a trivialization of the $H_{S\pr}$-torsor $G \times_{Q, x} S\pr \ra S\pr$ amounts to an element of $G(S\pr)$ lifting $x$, functorially in $S\pr$.

\item
For an $H$-torsor $T\pr$, a trivialization the associated $G$-torsor $T$ amounts to an identification of $T\pr$ with a sub-$H$-torsor of the trivial $G$-torsor, i.e., an identification of $T\pr$ with the fiber of $G \ra Q$ above an $S$-point of $Q$, functorially in $S$.

\item
For a $G$-torsor $T$, a trivialization of the $Q$-torsor $T/H$ amounts to a sub-$H$-torsor $T\pr \subset T$, i.e., an $H$-torsor $T\pr$ whose associated $G$-torsor identifies with $T$, functorially in $T$ and $S$.
 \qedhere
\eenum
\epf

\brem
The displayed squares remain Cartesian in a more general topos-theoretic setting, see \cite{IZ15}*{1.12, 1.16}.
\erem


\section{Lifting points of non-quasi-separated algebraic stacks} \lab{new-app-b}

The prevalent theme of this appendix is elimination of quasi-separatedness hypotheses from several known results. This is technical but useful: e.g.,~without our work here \Cref{open-supreme} would require various quasi-compactness and separatedness assumptions.

Our main goal is \Cref{LMB-revamp}, which extends \cite{Knu71}*{II.6.4} and \cite{LMB00}*{6.3} by allowing non-quasi-compact diagonals and by replacing fields by arbitrary Henselian local rings. The proof uses the same techniques as the above references and is modeled on the proof of \cite{LMB00}*{6.3}. After discussing several consequences of \Cref{LMB-revamp} (see, for instance, \Cref{cur-cor}), we combine them with \Cref{Groth-2} (which is the backbone of \S\ref{comparison-2}) to prove in \Cref{fppf-et,Gro68-revamp} the algebraic space analogues of the main results of \cite{Gro68c}*{Appendix}. The proofs in this part closely follow those of loc.~cit.

\bpp[$\SEC_n(\sX/\sY)$] \lab{SEC}
Let $S$ be a scheme, 
\[
\sX \ra \sY
\]
a representable separated morphism of $S$-algebraic stacks, and 
\[
(\sX/\sY)^n \ce \sX \times_\sY \dotsb \times_\sY \sX
\]
the $n$-fold fiber product for some $n \ge 0$. For an $S$-scheme $S\pr$, the groupoid $(\sX/\sY)^n(S\pr)$ is equivalent to that of $(n + 1)$-tuples $(s\pr, x_1, \dotsc, x_n)$ where $s\pr \colon S\pr \ra \sY$ is an $S$-morphism and $x_i$ is a section of the $S$-algebraic space morphism $S\pr \times_{s\pr, \sY} \sX \ra S\pr$. On $S\pr$-points, $(\sX/\sY)^n \ra \sY$ maps $(s\pr, x_1, \dotsc, x_n)$ to $s\pr$. Define the $S$-substack
\[
\SEC_n(\sX/\sY) \subset (\sX/\sY)^n
\]
by requiring the $x_i$ to be disjoint, i.e., 
\[
S\pr \times_{x_i, S\pr \times_{s\pr, \sY} \sX, x_j} S\pr = \emptyset \qq \text{for} \q i \neq j.
\] 
By \cite{SP}*{\href{http://stacks.math.columbia.edu/tag/03KP}{03KP}}, each $x_i$ is a closed immersion, so the inclusion $\SEC_n(\sX/\sY) \subset (\sX/\sY)^n$ is representable by open immersions. Thus, $\SEC_n(\sX/\sY)$ is an algebraic stack by \cite{SP}*{\href{http://stacks.math.columbia.edu/tag/06DC}{06DC}} (applied to the base change of a smooth atlas of $(\sX/\sY)^n$; compare also with \cite{LMB00}*{4.5 (ii)}). The formation of 
\[
\SEC_n(\sX/\sY) \subset (\sX/\sY)^n \ra \sY
\]
commutes with base change along any $\sY\pr \ra \sY$.
\epp

\bpp[$\ET_n(\sX/\sY)$] \lab{ET}
In the setup of \S\ref{SEC}, if $\sY$, and hence also $\SEC_n(\sX/\sY)$, is an algebraic space, then the symmetric group $\mathfrak{S}_n$ acts freely on $\SEC_n(\sX/\sY)$ by permuting the $x_i$. Sheafification of the constant presheaf $\fS_n$ to the constant group $(\fS_n)_S$ retains both the action and its freeness, so the quotient sheaf 
\[
\ET_n(\sX/\sY) \ce \SEC_n(\sX/\sY)/(\fS_n)_S
\]
is an $S$-algebraic space by \cite{SP}*{\href{http://stacks.math.columbia.edu/tag/06PH}{06PH}} (or, in many situations, already by \cite{LMB00}*{8.1.1 and 10.6}). By descent, if $S\pr$ is an $S$-scheme, then $\ET_n(\sX/\sY)(S\pr)$ is the set of pairs $(s\pr, x)$ where $s\pr\colon S\pr \ra \sY$ is an $S$-morphism and $x \subset S\pr \times_{s\pr, \sY} \sX$ is a closed subspace with $x \ra S\pr$ finite \'{e}tale of degree $n$.  

If $\sY$ is merely an algebraic stack, then we define the $S$-stack $\ET_n(\sX/\sY)$ by letting $\ET_n(\sX/\sY)(S\pr)$ be the groupoid of pairs $(s\pr, x)$ as in the previous sentence. The formation of 
\[
\SEC_n(\sX/\sY) \ra \ET_n(\sX/\sY) \ra \sY
\] 
commutes with base change along any $y\pr\colon \sY\pr \ra \sY$. We choose $y\pr$ to be a smooth cover by a scheme and combine the previous paragraph with \cite{SP}*{\href{http://stacks.math.columbia.edu/tag/05UL}{05UL}} (compare with \cite{LMB00}*{4.3.2}) to deduce that 
\begin{enumerate} [label={(\arabic*)}] 
\item
$\ET_n(\sX/\sY)$ is algebraic, 

\item \lab{SEC-ET-rep}
$\SEC_n(\sX/\sY) \ra \ET_n(\sX/\sY)$ and $\ET_n(\sX/\sY) \ra \sY$ are representable, and 

\item \lab{fet-deg-n}
$\SEC_n(\sX/\sY) \ra \ET_n(\sX/\sY)$ is finite \'{e}tale of degree $n!$. 
\eenum
Moreover, \ref{SEC-ET-rep}--\ref{fet-deg-n} imply the claims about $\ET_n(\sX/\sY) \ra \sY$ in 
\begin{enumerate} [label={(\arabic*)}] \addtocounter{enumi}{3}
\item \lab{ET-et}
If $\sX \ra \sY$ is smooth (resp.,~\'{e}tale), then so are 
\[
\SEC_n(\sX/\sY) \ra \sY \qq \text{and} \qq \ET_n(\sX/\sY) \ra \sY.
\]
\eenum
\epp

\blem \lab{ET-sch}
In the setup of \S\S\ref{SEC}--\ref{ET}, suppose that $\sX$ is a scheme and $\sY$ is an algebraic space.
\benum
\item \lab{ET-sch-a}
If $\sX \ra S$ is separated and $n > 0$, then the $S$-algebraic space $\ET_n(\sX/\sY)$ is separated.

\item \lab{ET-sch-b}
If $\sX \ra S$ is quasi-affine, $\sY \ra S$ is quasi-separated, $\sX \ra \sY$ is locally of finite type, and $n > 0$, then $\ET_n(\sX/\sY)$ is a scheme.
\eenum
\elem

\bpf \hfill
\benum
\item
By \cite{SP}*{\href{http://stacks.math.columbia.edu/tag/03HK}{03HK}}, $\Delta_{\sY/S}$ is separated, so 
\[
\qq(\sX/\sY)^n \ra S \qq \text{along with} \qq \SEC_n(\sX/\sY) \ra S
\]
inherits separatedness of $\sX \ra S$. Thus, \cite{EGAII}*{5.4.3~(i)} and \cite{EGAIV4}*{18.12.6} prove~that
\[
\qq\SEC_n(\sX/\sY) \times_S (\fS_n)_S \ra \SEC_n(\sX/\sY) \times_S \SEC_n(\sX/\sY), \quad\quad (c, \sigma) \mapsto (c, c\sigma)
\]
is a closed immersion, and the separatedness of $\SEC_n(\sX/\sY)/(\fS_n)_S \ra S$ follows.

\item
By \cite{SP}*{\href{http://stacks.math.columbia.edu/tag/03HK}{03HK}} and \cite{EGAIV4}*{18.12.12}, $\Delta_{\sY/S}$ is quasi-affine, so $(\sX/\sY)^n \ra S$ inherits quasi-affineness from $\sX \ra S$. Once the open immersion 
\[
\qq \SEC_n(\sX/\sY) \hra (\sX/\sY)^n
\]
is proved to be quasi-compact and hence also quasi-affine, $\SEC_n(\sX/\sY) \ra S$ will have to be quasi-affine, too, and \cite{SP}*{\href{http://stacks.math.columbia.edu/tag/07S7}{07S7}} (or \cite{SGA3Inew}*{4.1}) will give the claim.

The open $\SEC_n(\sX/\sY) \subset (\sX/\sY)^n$ is the complement of the union of the diagonals 
\[
\qq\Delta_{i, j}\colon (\sX/\sY)^{n - 1} \hra (\sX/\sY)^n \qq \text{for $i \neq j$.}
\]
Each $\Delta_{i, j}$ is a closed immersion that is of finite presentation by \cite{EGAIV1}*{1.4.3 (v)}. The inclusion of the complement of $\Delta_{i, j}$ is therefore a quasi-compact open immersion, and hence so is $\SEC_n(\sX/\sY) \hra (\sX/\sY)^n$. \qedhere
\eenum
\epf

\bpp[$\mathrm{P}_n(\sX/\sY)$] \lab{P}
In the setup of \S\S\ref{SEC}--\ref{ET}, suppose that $n > 0$ and $\sX$ is an algebraic space. As in \S\ref{ET}, $(\fS_n)_S$ acts (possibly non-freely) on the algebraic space $\SEC_n(\sX/\sY)$. We fix a closed embedding $\iota\colon (\fS_n)_S \hra (\GL_n)_S$, e.g., one furnished by permutation matrices, and set
\[
\mathrm{P}_n(\sX/\sY) \ce (\SEC_n(\sX/\sY) \times_S (\GL_n)_S)/(\fS_n)_S,
\]
where $(\fS_n)_S$ acts by 
\[
((c, g), \sigma) \mapsto (c \sigma, \iota(\sigma)\i g).
\]
By \cite{SP}*{\href{http://stacks.math.columbia.edu/tag/06PH}{06PH}} (or, in many situations, already by \cite{LMB00}*{8.1.1 and 10.6}), $\mathrm{P}_n(\sX/\sY)$ is an $S$-algebraic space (that depends on $\iota$). The formation of 
\[
\mathrm{P}_n(\sX/\sY) \ra \ET_n(\sX/\sY) \ra \sY
\]
commutes with base change along any representable $y\pr\colon \sY\pr \ra \sY$. Moreover, \S\ref{ET} proves
\begin{enumerate} [label={(\arabic*)}] 
\item \lab{Pn-GLn}
If $\sY$ is an algebraic space, then $\mathrm{P}_n(\sX/\sY)$ is a $\GL_n$-torsor over $\ET_n(\sX/\sY)$.
\eenum

We choose $y\pr$ to be a smooth surjection from an algebraic space and use \ref{Pn-GLn} to obtain

\begin{enumerate} [label={(\arabic*)}] \addtocounter{enumi}{1}
\item \lab{Pn-ETn-sm}
$\mathrm{P}_n(\sX/\sY) \ra \ET_n(\sX/\sY)$ is smooth.
\eenum
\epp

We are now ready for the main result of this appendix.

\bthm \lab{LMB-revamp}
Let $S$ be a scheme, $\sZ$ an $S$-algebraic stack, and $\Spec R$ an $S$-scheme with $R$ Henselian local. If the $S$-fiber of $\Delta_{\sZ/S}$ over the image of the closed point of $\Spec R$ is separated (in particular, if $\sZ$ is an algebraic space), then for every $z \in \sZ(R)$ there is a smooth (resp.,~\'{e}tale if $\sZ$ is an algebraic space) morphism 
\[
Z \ra \sZ
\]
with $Z$ an affine scheme such that $z$ lifts to $Z(R)$. 
\ethm

See \Cref{main-case} for an example of how the $S$-fibral hypothesis on $\Delta_{\sZ/S}$ is useful practice.

\bpf[Proof of \Cref{LMB-revamp} in the case when $\Delta_{\sZ/S}$ is separated and $R$ is a field]
We use \cite{SP}*{\href{http://stacks.math.columbia.edu/tag/04X5}{04X5}} (or \cite{Con07a}*{1.3}) to change\footnote{We later use the gained affineness of $S$ to ensure the quasi-affineness of $W \ra S$, which we need for \Cref{ET-sch}~\ref{ET-sch-b}. Quasi-separatedness of $S$ would suffice (such quasi-separatedness is an omnipresent convention in \cite{LMB00}).} $S$ to $\Spec \bZ$. The new diagonal $\Delta_{\sZ/\bZ}$ is separated, being the composite 
\[
\sZ \xra{\Delta_{\sZ/S}} \sZ \times_S \sZ \ra \sZ \times_{\bZ} \sZ
\] 
in which the second map inherits separatedness from $\Delta_{S/\bZ}$.

We choose a smooth (resp.,~\'{e}tale if $\sZ$ is an algebraic space) surjection 
\[
W \ra \sZ
\]
with $W$ a separated $S$-scheme and set 
\[
W_R \ce W \times_{\sZ, z} \Spec R.
\]
Since $W \ra S$ and $\Delta_{\sZ/S}$ are separated, so is $W \ra \sZ$. Thus, $W_R$ is a separated smooth nonempty $R$-algebraic space. We retain these properties of $W_R$ as well as smoothness (resp.,~\'{e}taleness) but not surjectivity of $W \ra \sZ$ by replacing $W$ by a suitable affine open $W\pr \subset W$ and $W_R$ by $W\pr \times_W W_R$; this gains affineness of $W \ra S$.

We apply \cite{EGAIV4}*{17.16.3~(ii)} to an \'{e}tale cover of $W_R$ by a scheme to find a finite separable field extension $R\pr/R$ and a $w\in W_R(R\pr)$. We set $n \ce [R\pr : R]$, and we use \cite{Knu71}*{II.6.2} to assume that 
$\Spec R\pr \xra{w} W_R$
is a monomorphism, and hence even a closed immersion by \cite{SP}*{\href{http://stacks.math.columbia.edu/tag/04NX}{04NX} and \href{http://stacks.math.columbia.edu/tag/05W8}{05W8}} (alternatively, by \cite{LMB00}*{A.2} and \cite{EGAIV4}*{18.12.6}). Then, by \S\ref{ET}, 
\[
w\in \ET_n(W_R/R)(R).
\]
Moreover, $w$ lifts to a $p \in \mathrm{P}_n(W_R/R)(R)$ thanks to the triviality of the $(\GL_n)_R$-torsor 
\[
\mathrm{P}_n(W_R/R) \times_{\ET_n(W_R/R), w} \Spec R.
\] 
The images of $w$ in $\ET_n(W/\sZ)(R)$ and of $p$ in $\mathrm{P}_n(W/\sZ)(R)$ lift $z$. As observed in \S\ref{ET} \ref{ET-et} and \S\ref{P} \ref{Pn-ETn-sm}, $\mathrm{P}_n(W/\sZ)$ inherits $\sZ$-smoothness from $W$, and $\ET_n(W/\sZ)$ also inherits $\sZ$-\'{e}taleness if $\sZ$ is an algebraic~space.

We use the argument above to pass to $\mathrm{P}_n(W/\sZ)$, and hence to reduce to the case when $\sZ$ is an algebraic space. Then we use the same argument and \Cref{ET-sch} \ref{ET-sch-a} to pass to $\ET_n(W/\sZ)$ and to assume further that $\sZ$ is separated. We repeat using \Cref{ET-sch} \ref{ET-sch-b} instead to assume that $\sZ$ is even a scheme, in which case an open affine $Z \subset \sZ$ through which $z$ factors suffices.
\epf

\emph{Proof of \Cref{LMB-revamp} in the case when $R$ is a field.}
If $\sZ$ is an algebraic space, then $\Delta_{\sZ/S}$ is separated as before. If not, then let 
\[
W \ra \sZ
\]
be a smooth surjection with $W$ a separated $S$-scheme. For $n \ge 1$, the $n$-fold fiber product $(W/\sZ)^n$ is a smooth $\sZ$-algebraic space with a $(\fS_n)_S$-action. By \cite{SP}*{\href{http://stacks.math.columbia.edu/tag/04TK}{04TK} and \href{http://stacks.math.columbia.edu/tag/04X0}{04X0}} (see also \cite{LMB00}*{4.3.1}), the quotient stack 
\[
(W/\sZ)^n/(\fS_n)_S
\]
is algebraic and 
\[
(W/\sZ)^n \ra (W/\sZ)^n/(\fS_n)_S
\] 
is a smooth surjection. Therefore, $(W/\sZ)^n/(\fS_n)_S$ inherits $\sZ$-smoothness from $(W/\sZ)^n$.
\bcl
The diagonal of $(W/\sZ)^n/(\fS_n)_S \ra S$ is separated.
\ecl

\bpf
We need to prove that 
\[
a\colon (W/\sZ)^n \times_S (\fS_n)_S \ra (W/\sZ)^n \times_S (W/\sZ)^n, \quad \quad \quad (x, \sigma) \mapsto (x, x\sigma)
\]
is separated. The diagonal $\Delta_{\mathrm{pr}_1 \circ a}$ of the composition of $a$ and 
\[
(W/\sZ)^n \times_S (W/\sZ)^n \xra{\mathrm{pr}_1} (W/\sZ)^n
\]
is
\[
(W/\sZ)^n \times_S (\fS_n)_S \xra{\Delta_a} F \xra{f} (W/\sZ)^n \times_S (\fS_n)_S \times_S (\fS_n)_S, \quad \text{where $f$ is a base change of $\Delta_{\mathrm{pr}_1}$}.
\]
Since 
\[
\Delta_{\mathrm{pr}_1 \circ a} = \id_{(W/\sZ)^n} \times \Delta_{(\fS_n)_S/S}
\]
and $\Delta_{\mathrm{pr}_1}$ is separated, $\Delta_a$ is a closed immersion.
\epf
To reduce to the previous case, it therefore suffices to lift $z$ to $((W/\sZ)^n/(\fS_n)_S)(R)$ for some $n \ge 1$. For this, since the formation of 
\[
(W/\sZ)^n \ra (W/\sZ)^n/(\fS_n)_S \ra \sZ
\]
commutes with base change along any $S\pr \ra S$, we choose $S\pr$ to be the point of $S$ below $z$ to reduce to the case when $\Delta_{\sZ/S}$ is separated. The proof of the previous case then lifts $z$ to $\ET_n(W/\sZ)(R)$ for some $n \ge 1$, and the presence of a $\sZ$-morphism 
\[
\ET_n(W/\sZ) \ra (W/\sZ)^n/(\fS_n)_S
\]
obtained from the construction finishes the proof. \QED

\emph{Proof of \Cref{LMB-revamp} in the general case.} The result follows by combining the proved field case with the following lemma.

\begin{sublem} \lab{lift-lift-2}
Let $S$ be a scheme, $\sZ$ an $S$-algebraic stack, and $\Spec R$ an $S$-scheme with $(R, \fm)$ Henselian local. For a smooth $S$-morphism 
\[
Z \ra \sZ
\]
with $Z$ an algebraic space, if the pullback $z_0 \in \sZ(R/\fm)$ of a $z \in \sZ(R)$ lifts to a $\wt{z}_0 \in Z(R/\fm)$, then $z$ lifts to a $\wt{z} \in Z(R)$ that pulls back to $\wt{z}_0$.
\end{sublem}

\emph{Proof.}
We set 
\[
Z_R \ce Z \times_{\sZ, z} \Spec R,
\]
so $\wt{z}_0$ gives rise to an element $\wt{z}_0 \in Z_R(R/\fm)$. To lift the latter to a desired section of 
\[
Z_R \ra \Spec R,
\]
we use the established field case of \Cref{LMB-revamp} to replace $Z_R$ by an \'{e}tale $Z_R$-scheme to which $\wt{z}_0$ lifts and then apply \cite{EGAIV4}*{18.5.17}.
\QEDD

\brem \lab{main-case}
We do not know if the $S$-fibral separatedness assumption on $\Delta_{\sZ/S}$ is necessary in \Cref{LMB-revamp}. Due to separatedness of group schemes over fields and \Cref{Delta-BG} \ref{Delta-BG-b}, this assumption is met if $\sZ = \bbB G$ for a flat and locally of finite presentation $S$-group algebraic space $G$ whose $S$-fiber over the image of the closed point of $\Spec R$ is a scheme.
\erem

\bcor \lab{lift-a-lot}
Let $S$ be a scheme, $\Spec R$ an $S$-scheme with $(R, \fm)$ Henselian local, and $\sZ$ an $S$-algebraic stack such that $\Delta_{\sZ/S}$ has a separated $S$-fiber over the image of the closed point of $\Spec R$. There is a separated $S$-scheme $Z$ and a smooth (resp.,~\'{e}tale if $\sZ$ is an algebraic space) $Z \ra \sZ$ for which 
\[
Z(R/\fa) \ra \sZ(R/\fa)
\]
is essentially surjective for every ideal $\fa \subset R$; if the number of isomorphism classes of objects of $\sZ(R/\fm)$ is finite, then $Z$ can be taken to be affine. 
\ecor

\bpf
We use \Cref{LMB-revamp} to build a $Z$ to which every $z\in \sZ(R/\fm)$ lifts and apply \Cref{lift-lift-2}.
\epf

\bcor \lab{hens-lem}
For a Henselian local ring $R$ and a smooth $R$-algebraic stack $\sZ$ such that $\Delta_{\sZ/R}$ has a separated $R$-fiber over the closed point of $\Spec R$, the pullback map 
\[
\sZ(R) \ra \sZ(R/\fa)
\]
is essentially surjective for every ideal $\fa \subsetneq R$. 
\ecor

\bpf
\Cref{lift-a-lot} reduces to the scheme case, which is a known variant of \cite{EGAIV4}*{18.5.17}.
\epf

We illustrate \Cref{hens-lem} with two special cases recorded in \Cref{lift-AVs,cur-cor}.

\bcor \lab{lift-AVs}
For a Henselian local ring $R$ and an ideal $\fa \subset R$, every principally polarized abelian scheme over $R/\fa$ arises as the base change of a principally polarized abelian scheme over $R$.
\ecor

\bpf
By \cite{FC90}*{(i) on p.~95}, for every $g \ge 0$ the moduli stack $\cA_g$ of principally polarized abelian schemes of relative dimension $g$ is separated and smooth over $\bZ$, so \Cref{hens-lem} applies to it.
\epf

The results of \ref{cur-cor}--\ref{Toe11-revamp} play a role in the proof of \Cref{sm-discrete}.

\bcor \lab{cur-cor}
For a Henselian local ring $R$, an ideal $\fa \subset R$, and a flat, locally of finite presentation $R$-group algebraic space $G$ whose special fiber is a scheme, 
\[
H^1(R, G) \ra H^1(R/\fa, G)
\] 
is surjective. 
\ecor

\bpf 
Thanks to \Cref{BG-sm} and \Cref{main-case}, \Cref{hens-lem} applies to $\bbB G$.
\epf

\brems
\remi \lab{sm-H1-inj}
For $\fa \subsetneq R$ as in \Cref{cur-cor} and $G$ a smooth $R$-group algebraic space, 
\[
H^1(R, G) \ra H^1(R/\fa, G)
\] 
is injective: fix $G$-torsors $T$ and $T\pr$ and apply \Cref{hens-lem}  to the fppf sheaf 
\[
\Isom_{G}(T, T\pr)\colon S \mapsto \{ \text{$G_S$-torsor isomorphisms $T_S \isomto T\pr_S$}\},
\]
which is a torsor under the inner form $\Aut_G T$ of $G$, and hence also a smooth $R$-algebraic space. Thus, if the special fiber of $G$ is a scheme, then 
$H^1(R, G) \ra H^1(R/\fa, G)$
is bijective (compare with \cite{SGA3IIInew}*{XXIV, 8.1~(iii)}). Likewise if $G$ is commutative: see \Cref{Gro68-revamp}. 

\remi
See \cite{Toe11}*{3.4} for an analogue of \Cref{cur-cor} with $R$ Henselian local excellent, $\fa$ the maximal ideal, and $G$ commutative, flat, and locally of finite presentation (but with arbitrary special fiber). Under these assumptions, loc.~cit.~also shows the bijectivity of the analogous pullback for $H^n$ with $n > 1$. The excellence assumption can be removed in many cases: 
\erems

\bprop\lab{Toe11-revamp}
For a Henselian local ring $(R, \fm)$, an ideal $\fa \subset \fm$, and a commutative flat $R$-group algebraic space $G$ of finite presentation, 
\[
H^1(R, G) \ra H^1(R/\fm, G)
\]
is surjective and 
\[
H^n(R, G) \ra H^n(R/\fa, G)
\]
is bijective for $n \ge 2$.
\eprop

\bpf
We may assume that $\fa = \fm$. We use \cite{EGAIV4}*{18.6.14 (ii)} to express $(R, \fm)$ as a filtered direct limit of Henselian local rings $(R_i, \fm_i)$ each of which is a Henselization of a finite type $\bZ$-algebra. By \cite{EGAIV2}*{7.8.3} and \cite{EGAIV4}*{18.7.6}, each $R_i$ is excellent. By \cite{SP}*{\href{http://stacks.math.columbia.edu/tag/07SK}{07SK} and \href{http://stacks.math.columbia.edu/tag/08K0}{08K0}}, there is an $i$ for which $G$ descends to a commutative flat $R_i$-group algebraic space $G_i$ of finite presentation. 

For each $j \ge i$, we set $G_j \ce (G_i)_{R_j}$. By the limit formalism for fppf cohomology, namely, by the analogue of \cite{SGA4II}*{VII, 5.9}, 
\[
\textstyle H^n(R, G) = \varinjlim_j H^n(R_j, G_j) \q \text{ and }\q H^n(R/\fm, G) = \varinjlim_j H^n(R_j/\fm_j, G_j).
\]
The claims therefore reduce to the case of an excellent $R$, which is the subject of \cite{Toe11}*{3.4}.
\epf

We now turn to the algebraic space analogues of the results of \cite{Gro68c}*{Appendix}. The proofs are analogous, too; we have decided to include them here because \Cref{Groth-2} is of major importance for \S\ref{comparison-2}, while the rest are quick consequences.

\blem \lab{SP-dig}
Let $S$ be a scheme, $j \colon T\pr \hra T$ a square-zero closed immersion of $S$-algebraic spaces, $G$ a commutative smooth $T$-group algebraic space, and $N$ the fppf sheaf $\Ker(G \ra j_*(G_{T\pr}))$. There is a quasi-coherent sheaf $\cF$ on $T$ such that 
\[
N(\wt{T}) = \Gamma(\cO_{\wt{T}}, \cF|_{\wt{T}})
\]
for every flat $T$-algebraic space $\wt{T}$.
\elem

\bpf
Let $\cC_{T\pr/T}$ be the conormal sheaf of $j$ (see \cite{SP}*{Definition \href{http://stacks.math.columbia.edu/tag/04CN}{04CN}}) and let $e\colon T \ra G$ be the unit section. Set 
\[
\cF \ce j_*\cH om_{\cO_{T\pr}}(j^*e^*\Omega_{G/T}, \cC_{T\pr/T}).
\] 
By \cite{SP}*{\href{http://stacks.math.columbia.edu/tag/04CQ}{04CQ}}, the formation of $\cC_{T\pr/T}$ commutes with base change along a flat $\wt{T} \ra T$, so \cite{SP}*{\href{http://stacks.math.columbia.edu/tag/061C}{061C} and \href{http://stacks.math.columbia.edu/tag/061D}{061D}} give the claim.
\epf

\blem \lab{Groth-2}
Let $S$ be a scheme, $T_0 \ra T$ a finite locally free morphism of $S$-algebraic spaces, and $G$ a commutative smooth $T$-group algebraic space. For $n \ge 0$, set 
\[
C^n \ce \Res_{T_n/T} (G_{T_n})
\]
where $T_n = T_0 \times_T \dotsc \times_T T_0$ with $n + 1$ factors of $T_0$, define the usual coboundary map $d^n\colon C^n \ra C^{n + 1}$, and set $Z^n \ce \Ker d^n$. The induced map
\[
d^n\colon C^n \ra Z^{n + 1}
\] 
is a smooth morphism of $T$-group algebraic spaces.  
\elem

\bpf
By \cite{SP}*{\href{http://stacks.math.columbia.edu/tag/05YF}{05YF} and \href{http://stacks.math.columbia.edu/tag/04AK}{04AK}} (see also \cite{Ols06}*{1.5}), $C^n$ and $Z^{n + 1}$ are locally of finite presentation $T$-group algebraic spaces. By \cite{SP}*{\href{http://stacks.math.columbia.edu/tag/04AM}{04AM}} (see also \cite{LMB00}*{4.15 (ii)}), it remains to argue that 
\[
d^n\colon C^n \ra Z^{n + 1}
\] 
is formally smooth. For this, we base change to assume that $T$ is affine, $j\colon T\pr \hra T$ is a closed subscheme defined by a square-zero ideal, and $a\pr \in C^n(T\pr)$ is such that $b\pr \ce d^n(a\pr) \in Z^{n + 1}(T\pr)$ lifts to a $b \in Z^{n + 1}(T)$. 

We need to lift $a\pr$ to an $a \in C^n(T)$ subject to $d^n(a) = b$. Formal smoothness of $C^n$ inherited from $G$ lifts $a\pr$ to an $\wt{a} \in C^n(T)$. We replace $b$ by $b - d^n(\wt{a})$ to reduce to the case when $a\pr = b\pr = 0$. We set 
\[
N \ce \Ker(G \ra j_*(G_{T\pr})),
\] 
so the cocycle $b$ is valued in the subsheaf $N \subset G$. 

By \Cref{SP-dig}, 
\[
H^{n + 1}(T_0/T, N) = H^{n + 1}(T_0/T, \cF)
\] 
for some quasi-coherent sheaf $\cF$ on $T_\fppf$, so $H^{n + 1}(T_0/T, N) = 0$ by \cite{Gro59}*{B, Lemme 1.1}. The existence of an $N$-valued cochain $a$ with $d^n(a) = b$ follows.
\epf

\blem \lab{pullb-cover}
Let $R$ be a Henselian local ring, $\fa \subsetneq R$ an ideal, $G$ a commutative smooth $R$-group algebraic space, $R_0$ an $R$-algebra that is finite free as an $R$-module, $T_0 = \Spec R_0$, and $T = \Spec R$. 
\benum
\item \lab{pullb-cover-b}
If $R$ is strictly Henselian and $n \ge 1$, then $H^n(T_0/T, G) = 0$.

\item \lab{pullb-cover-a}
The map 
\[
H^n(T_0/T, G) \ra H^n((T_0)_{R/\fa}/T_{R/\fa}, G)
\] 
is surjective for $n \ge 0$ and bijective for $n \ge 1$.
\eenum
\elem

\bpf
For an $x \in Z^n(R)$ with $n \ge 1$, we set 
\[
Z^n_x \ce \Spec R \times_{x, Z^n, d^{n - 1}} C^{n - 1} \qqqq \text{(so $Z^n_0 = Z^{n - 1}$).}
\]
By \Cref{Groth-2}, $Z^n_x$ is $R$-smooth, so \cite{EGAIV4}*{17.16.3 (ii)} proves \ref{pullb-cover-b} by showing that $Z^n_x(R) \neq \emptyset$. Also, \Cref{hens-lem} applied to $Z^n_x$ (resp.,~$Z^n_0$) proves the injectivity (resp.,~surjectivity) in \ref{pullb-cover-a}.
\epf

\bthm \lab{fppf-et}
Let $S$ be a scheme and $f\colon S_\fppf \ra S_\et$ the canonical morphism of sites. For a commutative smooth $S$-group algebraic space $G$, one has $\bbR^n f_* G = 0$ for every $n \ge 1$; in particular, 
\[
H^n_\et(S, G) = H^n_\fppf(S, G)
\] 
for such $G$ and $n$.
\ethm

\bpf
By \cite{Gro68c}*{11.1}, it suffices to prove that if $S = \Spec R$ with $R$ strictly Henselian local, then $H^n_\fppf(S, G) = 0$ for $n \ge 1$. For this, we combine \cite{Gro68c}*{11.2} with \Cref{pullb-cover} \ref{pullb-cover-b}.
\epf

\bthm \lab{Gro68-revamp}
Let $R$ be a ring, $\fa \subsetneq R$ an ideal, and $G$ a commutative smooth $R$-group algebraic space. If $\fa$ is nilpotent or if $R$ is Henselian local, then the pullback map
\[
H^n(R, G) \ra H^n(R/\fa, G)
\] 
is bijective for~$n \ge 1$.
\ethm

\bpf
For $j \colon \Spec R/\fa \hra \Spec R$, set $N \ce \Ker(G \ra j_*(G_{R/\fa}))$. Thanks to \Cref{hens-lem}, 
\[
0 \ra N \ra G \ra j_*(G_{R/\fa}) \ra 0
\]
is short exact in $(\Spec R)_\et$. Thus, \Cref{fppf-et} reduces us to proving that $H^n_\et(R, N) = 0$ for~$n \ge 1$. 

The case when $\fa$ is nilpotent reduces to case when $\fa^2 = 0$, in which $H^n_\et(R, N) = 0$ for $n \ge 1$ by \Cref{SP-dig} combined with the vanishing of quasi-coherent cohomology of affine schemes. 

In the case when $R$ is Henselian local, the analogue of \cite{Gro68c}*{11.2} for \'{e}tale cohomology and finite \'{e}tale covers reduces us to proving that $H^n(R_0/R, N) = 0$ for $n \ge 1$ and an \'{e}tale $R$-algebra $R_0$ that is finite as an $R$-module. We set $T_0 \ce \Spec R_0$ and $T \ce \Spec R$, consider sheaves of cochains for $T_0/T$ as in \Cref{Groth-2}, decorate them with subscripts to indicate the relevant groups, and note that by \Cref{hens-lem} 
\[
0 \ra C^n_N(R) \ra C^n_G(R) \ra C^n_{j_*(G_{R/\fa})}(R) \ra 0
\]
is short exact for $n \ge 0$. The resulting short exact sequence of complexes gives the exact sequence
\[
\dotsc  \ra H^{n}((T_0)_{R/\fa}/T_{R/\fa}, G) \ra H^{n + 1}(T_0/T, N) \ra H^{n + 1}(T_0/T, G) \ra H^{n + 1}((T_0)_{R/\fa}/T_{R/\fa}, G) \ra \dotsc,
\]
so the desired $H^n(T_0/T, N) = 0$ for $n \ge 1$ follows from \Cref{pullb-cover} \ref{pullb-cover-a}.
\epf


\end{document}